\documentclass[12pt,reqno]{amsart}
\usepackage{amsthm, amsmath,amssymb, diagrams}
\usepackage{color}
\numberwithin{equation}{section}
\usepackage{pgf,tikz}
\usetikzlibrary{arrows}
\usepackage{hyperref}

\theoremstyle{plain}

\newtheorem{thm}{Theorem}[section]
\newtheorem*{thmA}{Main Theorem}
\newtheorem{lem}[thm]{Lemma}
\newtheorem{prop}[thm]{Proposition}
\newtheorem{cor}[thm]{Corollary}

\newtheorem{conj}[thm]{Conjecture}

\theoremstyle{definition}

\newtheorem{defn}[thm]{Definition}
\newtheorem{exa}[thm]{Example}

\theoremstyle{remark}

\newtheorem{rmk}[thm]{Remark}

\newcommand{\g}[1]{\langle #1 \rangle}
\newcommand{\bp}[1]{\big{(} #1\big{)}}
\newcommand{\bt}[1]{
    \begin{tabular}[t]
    #1
    \end{tabular}
}
\newcommand{\ba}[1]{
    \begin{array}
    #1
    \end{array}
}
\newcommand{\bA}[1]{
    \begin{aligned}
    #1
    \end{aligned}
}
\newcommand{\bc}[1]{
    \begin{cases}
    #1
    \end{cases}
}
\newcommand{\bM}[1]{
    \begin{bmatrix}
    #1
    \end{bmatrix}
}
\newcommand{\bD}[1]{
    \begin{diagram}
    #1
    \end{diagram}
}
\newcommand{\ds}{\displaystyle}
\newcommand{\ts}{\textstyle}

\newcommand{\cd}{\cdot}
\newcommand{\cds}{\cdots}
\newcommand{\RR}{\mathbb{R}}

\newcommand{\CC}{\mathbb{C}}

\newcommand{\ZZ}{\mathbb{Z}}
\newcommand{\QQ}{\mathbb{Q}}

\newcommand{\KK}{\mathbb{K}}

\newcommand{\fa}{\text{ for all }}
\newcommand{\fsm}{~\text{ for some }~}
\newcommand{\otw}{\textup{otherwise}}
\newcommand{\tif}{\textup{if }}
\newcommand{\e}{\exists~}

\newcommand{\rw}{\rightarrow}

\newcommand{\Lrw}{\Leftrightarrow}

\newcommand{\uw}{\uparrow}
\newcommand{\ld}{\lambda}
\newcommand{\Ld}{\Lambda}
\newcommand{\af}{\alpha}
\newcommand{\Tt}{\Theta}
\newcommand{\fg}{\mathfrak{g}}
\newcommand{\fb}{\mathfrak{b}}
\newcommand{\fn}{\mathfrak{n}}
\newcommand{\fh}{\mathfrak{h}}

\newcommand{\mc}[1]{\mathcal{#1}}

\newcommand{\Rad}{\text{Rad}}

\newenvironment{enu}{\begin{enumerate}}{\end{enumerate}}
\newenvironment{eq}{\begin{equation}}{\end{equation}}

\begin{document}
\title[Weyl Modules over Affine Algebras]{
	On Weyl modules over affine Lie algebras in prime characteristic
}
\author{Chun-Ju Lai}
\address{
	Department of Mathematics, University of Virginia,
	Charlottesville, VA 22904
} 
\email{
	cl8ah@virginia.edu
}
\maketitle
\begin{abstract}
We construct a family of homomorphisms between Weyl modules for affine Lie algebras in characteristic $p$, 
which supports our conjecture on the strong linkage principle in this context. 
We also exhibit a large class of reducible Weyl modules beyond level one, for $p$ not necessarily small.

\end{abstract}
\section{Introduction}
\subsection{}
The theory of modular representations of reductive algebraic groups is well established \cite{RAG}. 
In this context, the Weyl modules play a fundamental role as the Verma modules for the category $\mathcal{O}$,
as Lusztig's conjecture provides a precise formula expressing the irreducible character in terms of characters of Weyl modules,
where the highest weights of these Weyl modules are controlled by the so-called strong linkage principle --
in other words, if a simple module is a composition factor of a Weyl module then the two highest weights are ``linked'' together -- this notion is introduced by Andersen \cite{And80}, following the notion of ``linkage principle'' given by Humphreys \cite{Hum71}.

Explicit homomorphisms between Weyl modules have been constructed by Carter-Lusztig \cite{CL74} and by Carter-Payne \cite{CP80} for type $A$.
For arbitrary types, homomorphisms have been constructed by Franklin \cite{Fr81, Fr88}. 
Franklin's work is based on Shapovalov's construction \cite{Sha72} on homomorphisms between Verma modules in characteristic 0.
He first constructed an integral version of Shapovalov elements, 
which gives rise to homomorphisms between the $\ZZ$-forms of Verma modules.
By modifying integral Shapovalov elements, one defines a map between the Weyl modules using reduction modulo $p$. 
On the one hand, a detailed calculation involving the contravariant form is needed to make sure the map is nonzero.
On the other hand, the BGG resolution is needed to ensure that the map is indeed a homomorphism.

We are interested in developing a theory of modular representations of affine Lie algebras.
So far there are very few results in the literature. 
For one, Mathieu \cite{Mat96} 
showed that the Steinberg modules are not irreducible in contrast with the finite types.
For another, DeConcini-Kac-Kazhdan \cite{DKK89} showed that  
the basic representations for untwisted affine ADE types are irreducible if and only if the characteristic $p$ of the underlying field
 does not divide the determinant of the Cartan matrix of the underlying finite root system.
Chari-Jing \cite{CJ01} further showed that the level one representations are irreducible if $p$ is coprime to the Coxeter number $h$.
For (possibly twisted) affine ADE types, Brundan-Kleshchev \cite{BK02} computed explicitly the determinant of the Shapovalov form on weight spaces of the basic representations.
When combined with \cite{DKK89}, it implies that the basic representation is irreducible if $p>3$ for (possibly twisted) affine DE types and if $p>h$ for (possibly twisted) affine A types.
The success of all these works above essentially hinges on concrete realizations of basic representations of affine Lie algebras via vertex operators. 
In contrast, little is known about the irreducibility of Weyl modules beyond level one.

A major difficulty in generalizing modular representation theory of finite types is 
that in an affine Weyl group there are no longest elements, which play an important role in the finite types.
For example, the strong linkage principle for modular finite types can be proved either by the Jantzen sum formula \cite{Jan77} under some assumption on $p$ or by Andersen's  approach \cite{And80} on the cohomology on line bundles.
In either way, the proof highly relies on the existence of the longest element of a Weyl group.

The main observation of this article is that Franklin's method can be generalized to the affine case 
as long as one establishes the affine version of the integral Shapovalov elements.
Here we state our main theorem. 
\begin{thmA}
Assume that $\ld,\mu$ are dominant integral weights such that $\mu = \ld - D \gamma$  for some positive real root $\gamma$ and positive integer $D$. There exists a nonzero homomorphism $V(\mu)\rw V(\ld)$  between Weyl modules of affine Lie algebras in prime characteristic if the following conditions hold.
\enu
\item $\mu$ is $\gamma$-mirrored to $\ld$  (cf.  \S\ref{pre:low}).
\item $\gamma$ is ``good'' (cf. \S \ref{defn:eta-good}).
\item $D$ is ``small enough'' (cf. Theorem \ref{thm:main}).
\endenu
\end{thmA}

In this article we also formulate a conjecture (Conjecture \ref{conj}) on the strong linkage principle for Weyl modules over affine Lie algebras.
To support our conjecture, we  construct  a large class of reducible Weyl modules beyond level one using the main theorem, for $p$ not necessarily small.
In particular, our results indicate that any clear formulation of a sufficient condition for the irreducibility of Weyl modules seems to be too subtle to formulate.

For type $\widetilde{A}_r$, we also give an upper bound for $p$ in terms of the level such that the main theorem applies. 
In particular, for type $\widetilde{A}_1$ and for the first few primes we give a table showing the lowest level of reducible Weyl modules that arise from our main theorem, 
and the complete list of quasi-simple weights (cf. \S \ref{defn:qs}) of level less than 150.

\subsection{}
The paper is organized as follows.
In Section 2 we introduce notations that are needed for the rest of this paper. 

In Section 3 we construct the integral Shapovalov elements for affine Lie algebras. 
These elements, by definition, give rise to nonzero homomorphisms between Verma modules for affine Lie algebras over $\ZZ$ if the two highest weights are linked by a reflection in the affine Weyl group. 

In Section 4 we show that if the two weights are $\gamma$-mirrored for some positive real root $\gamma$, 
the integral Shapovalov elements actually give rise to nonzero homomorphisms between Verma modules in the characteristic $p$, under a condition on the distance between the two highest weights. 
We then show that condition can be weakened by modifying the integral Shapovalov elements, provided $\gamma$ is ``good''.

In Section 5 we show that if two $\gamma$-mirrored weights satisfy an additional condition on their distance, 
an analogous nonzero homomorphisms between the corresponding Weyl modules can be constructed. 
The contravariant forms and the BGG resolution for Kac-Moody algebras are used in the proof. 
Another essential tool here is Lemma \ref{lem:CuZ}, which is a variant of the Shapovalov determinant formula. 

In Section 6  we relate our main theorem to our conjecture on the strong linkage principle.
We also describe the candidates for reducible Weyl modules. 
Using our main theorem, we show that most candidates are reducible except for some rare ``quasi-simple weights''. 
For type $\widetilde{A}_1$ we computed for the first few primes a list of reducible Weyl modules that can be detected by our main Theorem;
for type $\widetilde{A}_r$ we also give an upper bound for $p$ in terms of the level such that our main theorem applies.

\subsection*{Acknowledgments} 
I am grateful to my advisor, Professor Weiqiang Wang, for his guidance and support. 
I would like to thank Professor Leonard Scott for insightful discussions and bringing Franklin's paper to our attention.
I must thank Professor James Franklin for providing his unpublished thesis paper.
I  thank Professor Toshiyuki Tanisaki for pointing out an ambiguity in the paper.
 I would also like to thank Professor Jim Humphreys and Shrawan Kumar for useful comments.

This research is partially supported by the GRA fellowship of Wang's NSF grant DMS-1101268. 
We thank Institute of Mathematics, Academia Sinica, Taipei for providing an excellent
working environment and support, where part of this project was carried out.
\section{Preliminaries}

	\subsection{Affine Kac-Moody Algebras}
Let $\fg_\CC=\fg_\CC(A)$ be the affine Kac-Moody algebra over $\CC$ associated to a symmetrizable generalized Cartan matrix (GCM) $A= (a_{ij})_{i,j \in I}$ of affine type and of rank $r$ where $I=\{0,1,\ldots,r\}$. Following  \cite{Kac90} and \cite{Car05},
let $\vec{a}= (a_i)_{i\in I}$ and $\vec{c} = (c_i)_{i\in I}$ be the minimal positive integer vectors  such that
$A\vec{a} = 0 = \vec{c}^{~t} A$.
Let  $\left(\fh_\CC,(\af_i)_{i\in I},(h_i)_{i\in I}\right)$ be a minimal realization. 
Here $\fh_\CC$ is the $\CC$-span of $h_0,\ldots, h_r$, and $d \in \fh_\CC$ is a fixed scaling element such that $\af_i(d) = \delta_{i,0}$.
Let  $\fn_\CC^+$ (and $\fn_\CC^-$) be the subalgebra of $\fg_\CC$ generated by $(e_i)_{i\in I}$ (and $(f_i)_{i\in I}$, respectively). 
Let $\fb_\CC^\pm = \fh_\CC\oplus \fn_\CC^\pm$ be the corresponding Borel subalgebras.
Let $(,):\fg_\CC\times\fg_\CC\rw\CC$ be a non-degenerate invariant symmetric bilinear form which induces a non-degenerate bilinear form on $\fh_\CC^*$. 

Let $\delta = \sum_{i\in I} a_i \af_i$ be the basic imaginary root, 
let $c= \sum_{i\in I} c_i h_i$ be the central element, and fix a vector $\rho\in\fh_\CC^*$ satisfying 
$\rho(h_i) = 1 \fa i \in I$.
Let $h = \sum_{i\in I} a_i$ be the  Coxeter number, and let $h^\vee = \sum_{i\in I} c_i$ be the dual Coxeter number.
Let $X=\{\ld \in \fh_\CC^*\mid \ld(h_i) \in \ZZ \fa i \in I\}$ be the set of integral weights, 
and let $X^+=\{\ld \in X\mid  \ld(h_i) \geq 0 \fa i \in I\}$ be the set of dominant integral weights. 
Let $\varpi_0, \ldots, \varpi_r\in \fh_\CC^*$ be the fundamental weights, 
that is, each $\varpi_j$ is an element in $\fh_\CC^*$ such that  
$\varpi_j(h_i) = \delta_{i,j}$ and $\varpi_j(d)=0$. 
Therefore 
\[
X^+ = \left\{ \ld = \sum_{i=0}^r \xi_i \varpi_i + \xi \delta \in \fh_\CC^* ~\middle|~ \xi_i \in \ZZ_{\geq0} \fa i \in I,\xi\in\CC \right\}.
\]


	\subsection{Affine root systems }\label{pre:root}
Let $\Phi, \Phi^+, \Phi_\text{re}, \Phi^+_\text{re}, \Phi_\text{im}$ and $\Phi^+_\text{im}$ be the set of roots, positive roots, real roots, positive real roots, imaginary roots and positive imaginary roots of the affine Kac-Moody algebra $\fg_\CC$, respectively. 
Let $\Phi_0$ be the root system generated by $\{\af_i\}_{i=1}^r$.
Let $Q=\sum_{i\in I} \ZZ  \af_i$ (and $Q^+ =\sum_{i\in I} \ZZ_{\geq0} \af_i$) be the set of (non-negative) integer linear combinations of roots of $\fg_\CC$. 
For each element $\gamma = \sum_{i\in I} g_i \af_i \in Q$ for some $g_i \in \ZZ$, we set 
\eq\label{eq:hr}
h_\gamma= \sum_i g_i^\vee h_i, ~\text{ where }~ g_i^\vee = (a_i/c_i) g_i,
\endeq 
and let ht$\gamma = \sum_{i\in I} g_i$ be the height of $\gamma$. A simple root $\af_\eta$ is said to be \textit{occurring} in $\gamma$ if $g_\eta \neq 0$.

Let $\fh_\RR$ be the $\RR$-span of $\{h_i\}_{i\in I}$ and $d$. 
To each real root $\gamma \in \Phi_{\text{re}}$ we assign a reflection $s_\gamma : \fh_\RR^* \rw \fh_\RR^*$ given by
\[
s_\gamma(\ld) = \ld-\g{\ld,\gamma^\vee}\gamma 
~\text{ where }~
\g{\ld,\gamma^\vee}= \frac{2(\ld,\gamma)}{(\gamma,\gamma)}.
\]
In addition, for each $m\in\ZZ$ we assign an affine reflection $s_{\gamma,m}:\fh_\RR^*\rw\fh_\RR^*$ given by
\[
s_{\gamma,m}(\ld) = s_\gamma(\ld)+m\gamma.
\]
The affine reflection $s_{\gamma,m}$ is the reflection with respect to the hyperplane 
\[
H_{\gamma,m} = \{\ld \in \fh_\RR^*\mid  \g{\ld+\rho,\gamma^\vee} = m\}.
\]
	\subsection{Affine Weyl groups} \label{pre:Weyl grp}
Let $W$ be the Weyl group of $\fg_\CC$ generated by $s_\af$ with $\af \in \Phi_{\text{re}}$. 
It is known that $W$ is also generated by the affine reflections $s_{\af,m}$ with $\af \in \Phi_0, m\in \ZZ$.
We introduce the dot action of $W$ on $\fh_\RR^*$ by
\[
w\cd\ld = w(\ld+\rho)-\rho.
\]
In particular, for each real root $\af$ we have
\[
s_{\af,m}\cd\ld = \ld - (\g{\ld+\rho,\af^\vee}-m)\af.
\]

Let $\theta$ be the highest root in $\Phi_0$. Then $(W,\mathcal{S})$ is a Coxeter system where the generating set $\mathcal{S}=\{s_{\theta,1},s_{\af_1},\ldots,s_{\af_r}\}$.
A reduced expression for $w\in W$ is a product $w=t_1\cds t_N$ of elements $\{t_i\}_{i=1}^N$ in $\mathcal{S}$ such that the number $N$ is minimal among all such expressions for $w$. Here $l(w)= N$ is called the length of $w$. 

	\subsection{Verma modules}\label{pre:Verma}
For any Lie algebra $\mathfrak{a}$ let $U(\mathfrak{a})$ be its universal enveloping algebra. 
Let $U = U(\fg_\CC)$,
and let  $\fh_\ZZ$ be the $\ZZ$-span of $\{ d, h_i ~|~ i\in I\}$. Denote the $\ZZ$-form of $U$ by
\[
U_\ZZ = \left\langle e_i^{(n)}, f_i^{(n)}, {h \choose n} ~\middle|~ i\in I, n\in \ZZ_{>0}, h\in \fh_\ZZ \right\rangle,
\]
where $e_i^{(n)}= e_i^n/n!, f_i^{(n)}=f_i^n/n!$  and ${h \choose n}= h(h-1)\cds(h-n+1)/n!$.
We denote by $U_\ZZ^-$ and $U_\ZZ^+$ the subrings of $U$ generated by $f_i^{(n)}$ and by $e_i^{(n)}$, respectively.
Let $\fg_\ZZ= \fg_\CC \cap U_\ZZ, \fn^\pm_\ZZ=\fn_\CC^\pm \cap U^\pm_\ZZ$ and $\fb^\pm_\ZZ= \fn_\ZZ^\pm \oplus \fh_\ZZ$.


Let $\KK$ be an algebraically closed field, 
let $\fg_\KK, \fn^\pm_\KK, \fh_\KK, U_\KK$ and $U^\pm_\KK$ be the tensor product of $\KK$ with $\fg_\ZZ, \fn^\pm_\ZZ, \fh_\ZZ, U_\ZZ$ and $U^\pm_\ZZ$, respectively. 
For each $\ld \in \fh^*_\ZZ$ we denote the Verma module (as a $U_\KK$-module) by 
\[
M(\ld)_\KK = U_\KK/I(\ld)_\KK,
\]
where $I(\ld)_\KK = I(\ld)_\ZZ \otimes \KK$ and 
\eq\label{eq:IZ}
I(\ld)_\ZZ = \sum_{j\in I,m\geq 1} U_\ZZ e_j^{(m)} +\sum_{j\in I,m\geq 1} U(\fb^-_\ZZ) { h_j - \ld(h_j) \choose m}.
\endeq
Let $v^+_{\ld,\KK}=1+I(\ld)_\KK$ be a highest weight vector of $M(\ld)_\KK$
 and let $v^+_\ld = v^+_{\ld,\CC}$.
The Verma module $M(\ld)_\CC$ contains a unique maximal submodule $N(\ld)_\CC$ and a unique simple quotient $L(\ld)_\CC = M(\ld)_\CC/N(\ld)_\CC$.

In \cite[Theorem 2]{KK79} Kac and Kazhdan proved the strong linkage principle for symmetrisable Kac-Moody algebras.
\prop\label{prop:lp}
$L(\mu)_\CC$ is a composition factor of $M(\ld)_\CC$ if and only if $\mu =\ld$ or $\mu = \mu_0 \uw \cds \uw \mu_N = \ld$ where
\eq\label{SL:aff ord}
x \uw y \Lrw \text{ there exists } \bc{n\in \ZZ_{>0},\\ \beta \in \Phi^+} \text{such that } \bc{y-x= n\beta,\\ n(\beta,\beta) = 2(y+\rho, \beta).}
\endeq
\endprop
\rmk
The notion of ``linkage principle'' --  all composition factors of an indecomposable have linked highest weights -- is introduced by Humphreys in \cite{Hum71}. The notion of the strong linkage principle -- all composition factors have highest weights strongly linked to a certain weight -- was then defined by Andersen \cite{And80},  though Jantzen had earlier worked out some of the ideas in his own notation.
\endrmk
	\subsection{Integral bases for $U_\ZZ$} \label{pre:PBW}
Based on Garland's work \cite{Gar78} on the loop algebras, 
Mitzman constructed a PBW-type basis for all (possibly twisted) affine Lie algebras.
Here we first recall the root vectors defined in \cite{Mit85};
define such a basis for each choice of root vector; then fix a choice (cf. Proposition \ref{prop:scale}) that significantly simplifies the computation. 

Assume that $\fg_\CC$ is of type $X_r^{(k)}$. Let $\sigma$ be the Dynkin diagram automorphism of order $k$ which fixes vertex 0.
For $1\leq i \leq r$ and $n\neq 0$, we set
\[
X_{n\delta,i} = 
\bc{
\frac{1}{2} X_{n\delta,\af_i} &\textup{if } (\af_i, \sigma(\af_i)) \textup{ is odd and } n \textup{ is even},
\\
X_{n\delta,\af_i} &\textup{otherwise}.
}
\]
Here $X_{n\delta,\af_i}$ is the vector defined in \cite[(3.9.1), (3.9.6)]{Mit85} (was denoted by $X_{n\imath,\textup{a}}$).
We then choose $X_\gamma$ for $\gamma \in \Phi_{\text{re}}^+$ such that the set 
$
\{ X_\gamma ~|~ \gamma \in \Phi_{\text{re}}^+\} 
\cup 
\{ X_{n\delta,i} ~|~1\leq i \leq r, n\neq 0\}
\cup
\{
h_0 , h_1, \ldots, h_r, d
\}$ is a Chevally basis of $\fg_\CC$ in the sense of \cite[Definition 2.2.25]{Mit85}.

Fix a total ordering $\leq$ on the set 
$\Phi^+_M = \{(\gamma, i) \in \Phi^+ \times \ZZ ~|~ 1\leq i \leq \dim  \fg_{\CC,\gamma} \}$. 
For now we fix root vectors 
\eq\label{eq:roots}
\{f_{\gamma,i},~  e_{\gamma,i} ~|~ (\gamma,  i) \in  \Phi_M^+ \}, 
\endeq
where
$e_{\gamma,1} = X_{\gamma}$,
$f_{\gamma,1} = X_{-\gamma}$
if $\gamma \in \Phi^+_{\textup{re}}$;
$e_{\gamma,i} = X_{\gamma,i}$,
$f_{\gamma,i} = X_{-\gamma,i}$
if $\gamma \in \Phi^+_{\textup{im}}$. 
For short, we write $f_\gamma = f_{\gamma,1}$ and $e_\gamma = e_{\gamma,1}$ for $\gamma \in  \Phi^+_\text{re}$;
write $e_i = e_{\af_i}$ and $f_i = f_{\af_i}$ for $i\in I$.
Define polynomials $\Ld_s = \Ld_s(x_1, x_2, x_3, \ldots)$ over $\CC[x_1, x_2, x_3, \ldots][[\zeta]]$ for $s\geq0$ by 
\[
\sum_{s\geq 0}\Ld_s \zeta^s = \textup{exp}\biggl(\sum_{j\geq 1} x_j \frac{\zeta^{j}}{j}\biggr).
\]
For $(\gamma,i) \in \Phi_M^+, h \in \fh_\ZZ, s\geq 0$, by a slight abuse of notation, set
\[
f_{\gamma,i}^{(s)}
=
\Ld_s(f_{\gamma,i}, f_{2\gamma,i}, f_{3\gamma,i}, \ldots),
e_{\gamma,i}^{(s)}
=
\Ld_s(e_{\gamma,i}, e_{2\gamma,i}, e_{3\gamma,i}, \ldots),
\textup{ and }
h^{(s)} = \Ld_s(h,h, \ldots).
\]
Note that by \cite[(4.2.9 -- 10)]{Mit85}, $f_{\gamma,1}^{(s)} = f_{\gamma}^{s}/s!$ and $e_{\gamma,1}^{(s)} = e_{\gamma}^{s}/s!$ are the usual divided powers if $\gamma \in \Phi_{\text{re}}^+$, and $h^{(s)} = {h+s-1 \choose s}$ for $h \in \fh_\ZZ$.
Let $\mathcal{P}$ be the set of functions mapping from $\Phi_M^+$ to $\ZZ_{\geq0}$ with finite support. 
For each $\pi\in\mathcal{P}$ with finite support $\{(\gamma_1,i_1) < \cds < (\gamma_m,i_m)\}$ we denote the degree of $\pi$ by
$\deg\pi = \sum_{j=1}^m \pi(\gamma_j, i_j)\in \ZZ_{\geq 0}$, the weight of $\pi$ by 
wt$(\pi) = \sum_{j=1}^N \pi(\gamma_j, i_j)\gamma_j\in Q^+$. We define 
\[
\ba{{ll}
f_\pi &= f_{\gamma_1,i_1}^{(\pi(\gamma_1, i_1))} \cds f_{\gamma_m,i_m}^{(\pi(\gamma_m, i_m))} \in U^-_\ZZ,\textup{ and } \\
e_\pi &= e_{\gamma_1,i_1}^{(\pi(\gamma_1,i_1))} \cds e_{\gamma_m,i_m}^{(\pi(\gamma_m, i_m))} \in U^+_\ZZ.
}
\]
Similarly, for each tuple $\varphi = (\varphi_0, \varphi_1, \ldots \varphi_{r}, \varphi_{r+1}) \in \ZZ_{\geq0}^{r+2}$ we assign
\[
h_\varphi = h_0^{(\varphi_0)}\cds h_r^{(\varphi_r)} d^{(\varphi_{r+1})} \in U(\fh_\ZZ).
\]
\begin{lem}(Mitzman)
$U_\ZZ$ has a Garland-Kostant-type $\ZZ$-basis 
\[
\{f_\pi h_\varphi e_{\pi'}\mid  \pi,\pi' \in \mathcal{P}, \varphi \in \ZZ_{\geq0}^{r+2}\}.
\]
As a consequence, $U$ has a PBW-type $\CC$-basis 
\[
\{F_\pi H_\varphi E_{\pi'}\mid  \pi,\pi' \in \mathcal{P}, \varphi \in \ZZ_{\geq0}^{r+2}\},
\]
where 
\[
F_\pi = f_\pi \prod_{j=1}^m \pi(\gamma_j, i_j)!,
\quad
E_\pi = e_\pi \prod_{j=1}^m \pi(\gamma_j, i_j)!,
\quad
H_\varphi = h_\varphi \prod_{j=0}^{r+1} \varphi_j!.
\]
\end{lem}
\proof
By \cite[Theorem 4.2.6]{Mit85}, the set of  monomials $f_\pi h_\varphi e_{\pi'}$ forms a $\ZZ$-basis of the subalgebra of $U$ generated by $f_{\gamma,i}^{(s)}, e_{\gamma,i}^{(s)}, h_j^{(s)}$ and $d^{(s)}$ for $(\gamma,i) \in \Phi_M^+, 0 \leq j \leq r, s\geq0$. See \cite[Remark 4.2.7, Corollary 4.4.12]{Mit85} for the identification of this subalgebra with the algebra $U_\ZZ$ defined in Section \ref{pre:Verma}.
\endproof
We are able to express each $u \in U(\fb_\CC^-)$ as
\[
u = \sum_{\pi\in\mathcal{P}} F_\pi Q_\pi,
\] 
where $Q_\pi \in U(\fh_\CC) \simeq \CC[h_0, \ldots, h_{r},d]$. 
We further denote the  evaluation of $u$ at $\ld \in \fh_\CC^*$ by
\[
u(\ld) = \sum_{\pi\in\mathcal{P}} F_\pi Q_\pi(\ld) \in U(\fn_\CC^-).
\] 
Now we explain how we pick a choice of the root vectors and hence an integral basis of $U_\ZZ$.
For $\omega\in\mathcal{P}, q\in\ZZ_{>0}$ and a simple root $\af$, we may rewrite $f_\af^q F_\omega$,
 using relations in $U$, as a $\CC$-linear combination of monomials of the form
\[
F_\pi f_\af^{q-i}, ~\textup{ where }~ \pi(\af,1) = 0, i\in\ZZ.
\]
That is, we ``move'' $f_\af$ to the very right end. This operation shows up frequently when it comes to solving the Shapovalov-type equations in Section 3. 
Now we define $C_{\omega,\pi,i, q} \in\CC$ be the coefficients such that
\[
f_\af^m F_\omega = \sum_{\substack{\pi \in \mathcal{P}\\ \pi(\af,1)=0}}\sum_{i\in\ZZ} F_\pi f_\af^{m-i} C_{\omega,\pi,i,m}.
\]
\prop\label{prop:scale}
We can rescale the root vectors (cf. \eqref{eq:roots}) so that there are
polynomials $C_{\omega,\pi,i} \in \ZZ[m]$ satisfying 
$C_{\omega,\pi,i,m} = C_{\omega,\pi,i}(m)$ for all $\omega, \pi, i$ and $m$.
\endprop
\proof
The proof is exactly the same as in the finite case (see \cite[Lemma 4.1]{Fr88}). 
\endproof
From now on, by a little abuse of notation, fix $\{f_{\gamma,i},~  e_{\gamma,i} \} $ to be a set of root vectors as in Proposition \ref{prop:scale}, and hence we have
\[
f_\af^m F_\omega = \sum_{\substack{\pi \in \mathcal{P}\\ \pi(\af)=0}}\sum_{i\in\ZZ} F_\pi f_\af^{m-i} C_{\omega,\pi,i}(m)
\]
for each $m\in \ZZ_{>0}, \omega \in \mathcal{P}$ and simple root $\af$.  
\exa \label{exaC}
Consider the affine Kac-Moody algebra of type $\widetilde{A}_1$ with GCM 
\[
A=\bM{2&-2\\-2&2}.
\] 
Let $\gamma = \af_0+\delta$ and pick root vectors $f_0, f_1, f_{\delta,1}$ and $f_\gamma$ satisfying
\[
f_{\delta,1} = f_0f_1 - f_1f_0 \textup{ and } f_\gamma = f_0 f_{\delta,1} - f_{\delta,1} f_0.
\]
Fixing $\af = \af_0$ and $F_\omega = f_0^2f_1 = f_1 f_0^2  + 2 f_{\delta,1} f_0 + f_\gamma$, we have
\[
f_0^m F_\omega = ( f_1 f_0^2  + (m+2)  f_{\delta,1}f_0 + \textstyle{{m+2 \choose 2}} f_\gamma) f_0^m,
\]
that is, 
\[
C_{\omega,\pi, i} (m)
= \bc{
1 &\tif F_\pi = f_1, i=-2, 
\\
m+2 &\tif F_\pi = f_{\delta,1}, i=-1, 
\\
{m+2 \choose 2} &\tif F_\pi = f_\gamma, i=0,
\\
0&\otw.
}
\]
Note that ${m+2 \choose 2} \not\in \ZZ[m]$, so we have to choose again $f_\gamma = \frac{1}{2}(f_0 f_{\delta,1} - f_{\delta,1} f_0)$ to make $C_{\omega,\pi, i} (m) \in \ZZ[m]$ as mentioned in Proposition \ref{prop:scale}.
\endexa
	\subsection{Weyl modules}\label{pre:Weyl}
Note that in the finite dimensional framework there are two constructions of the Weyl modules -- one by reduction modulo $p$ from known modules in characteristic $0$, and one by taking dual of some induced modules. In this paper we adopt the definition of Weyl modules by reduction modulo $p$.
From now on let $\KK$ be an algebraically closed field of characteristic $p>0$. 
For each $\ld \in X^+$, recall that $M(\ld)_\CC = U v_{\ld}^+$ so that we have 
\[
L(\ld)_\CC = U \overline{v}_{\ld}^+ ~\text{ where }~ \overline{v}_{\ld}^+= v_{\ld}^+ + N(\ld)_\CC.
\] 
Then we define its $\ZZ$-form $L(\ld)_\ZZ= U_\ZZ\overline{v}_{\ld}^+$ and define the \textit{Weyl module} 
\[
V(\ld) = L(\ld)_\ZZ \otimes \KK. 
\]
	\subsection{Nearest lower reflections}\label{pre:low}
Fix a positive real root $\gamma\in \Phi_{\text{re}}^+$ and a positive integer $e$, 
for each $\ld \in X$ we set $D' = \g{\ld+\rho,\gamma^\vee}\in\ZZ$. 
There are uniquely determined integers $M,D \in \ZZ$ such that $D' = Mp^e + D$ such that $0 \leq D < p^e$.
Set $R_{\gamma,e}= s_{\gamma,Mp^e} $, hence
\[
R_{\gamma,e}\cd\ld = s_{\gamma,Mp^e}\cd\ld = \ld - D\gamma.
\]
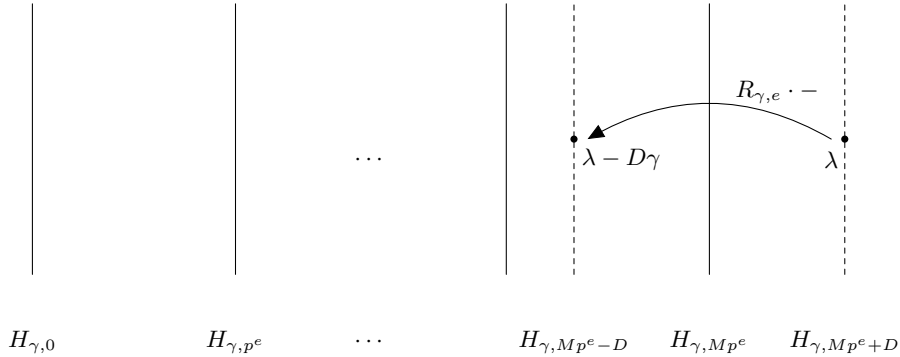
\[
\begin{tikzpicture}[line cap=round,line join=round,>=triangle 45,scale = 0.9]
\draw [-] (0,1) -- (0,5);
\draw [-] (3,1) -- (3,5);
\draw [-] (7,1) -- (7,5);
\draw [-,dash pattern=on 2pt off 2pt] (8,1) -- (8,5);
\draw [-] (10,1) -- (10,5);
\draw [-,dash pattern=on 2pt off 2pt] (12,1) -- (12,5);
\draw [-] (13,1) -- (13,5);
\path [->] (11.8,3)  edge[bend right]   (8.2,3);
\begin{scriptsize}
\draw[color=black] (0,0) node {$H_{\gamma,0}$};
\draw[color=black] (3,0) node {$H_{\gamma,p^e}$};
\draw[color=black] (8,0) node {$H_{\gamma,Mp^e-D}$};
\draw[color=black] (10,0) node {$H_{\gamma,Mp^e}$};
\fill [color=black] (12,3) circle (1.5pt);
\draw[color=black] (11.8,2.7) node {$\lambda$};
\fill [color=black] (8,3) circle (1.5pt);
\draw[color=black] (8.7,2.7) node {$\lambda-D\gamma$};
\draw[color=black] (5,0) node {$\cdots$};
\draw[color=black] (5,2.7) node {$\cdots$};
\draw[color=black] (12,0) node {$H_{\gamma,Mp^e+D}$};
\draw[color=black] (11,3.7) node {$R_{\gamma,e}\cdot -$};
\end{scriptsize}
\end{tikzpicture}
\]
In words, $R_{\gamma,e}$ is the \textit{nearest lower $p^e$-reflection}. That is, among those $p^e$-hyperplanes of the form $H_{\gamma, Mp^e}$ with $M\in\ZZ$, $R_{\gamma,e}$ is the reflection with respect to the hyperplane that is nearest to $\ld$.

Let $\ld,\mu \in X^+$ be two dominant integral weights. 
$\mu$ is said to be $(\gamma,e)$-\textit{mirrored} to $\ld$ if $\mu = R_{\gamma,e}\cd\ld$ and we write $\mu\uw_\gamma^e \ld$.
Finally, $\mu$ is said to be $\gamma$-\textit{mirrored} to $\ld$ if $\mu$ is $(\gamma,e)$-\textit{mirrored} to $\ld$ for some positive integer $e$ and we write $\mu \uw_\gamma \ld$.

\rmk
In modular finite case, the blocks were classified by Donkin \cite{Don80} using the theory of $G_n T$-modules. 
Here $G_n$ is the $n$th Frobenius kernel in the algebraic group $G$ and $T$ is a maximal torus in $G$.
In the representation theory for $G_nT$, the building blocks of strong linkage are given by the $\gamma$-mirrored condition. 
The strong linkage principle in this context is proved by Doty \cite{Dot89} and he showed that it implies the strong linkage principle for the modular finite case.  
This is why our seemingly strange $\gamma$-mirrored condition is the reasonable format.
\endrmk
\section{Construction of Integral Shapovalov elements}  \label{sec:ShElt}
In this section we construct an integral version of Shapovalov elements for the affine case. The argument is adapted from Franklin's construction \cite{Fr88} for the finite case. The idea is to study closely the coefficients in every Shapovalov equation in terms of the PBW-type basis for affine Lie algebras. 
\subsection{Generic Shapovalov Elements}
As in \cite{KK79}, it is easy to verify that if $\chi \in H_{\gamma,D}$ then any nonzero homomorphism $M(\chi-D\gamma)_\CC \rw M(\chi)_\CC$ is an embedding and is of the form 
$v^+_{\chi-D\gamma} \mapsto uv^+_\chi \fsm u\in U(\fb_\CC^-)_{-D\gamma}.$
On the other hand,  for a fixed positive real root $\gamma$ and a fixed positive integer $D$, any generic Shapovalov element $S$ (with respect to $(\gamma,D)$) gives rise to nonzero homomorphisms $M(\chi-D\gamma)_\CC \rw M(\chi)_\CC$ by $v^+_{\chi-D\gamma} \mapsto Sv^+_\chi$ if $\chi \in H_{\gamma,D}$.
\defn
An element $S\in U(\fb_\CC^-)$ is called a \textit{(generic) Shapovalov element} 
(with respect to a positive real root $\gamma$ and a positive integer $D$) 
if the conditions (S$1$)-(S$3$) below hold:
\enu
\item[(S$1$)] $S \in U(\fb_\CC^-)_{-D\gamma}$; 
\item[(S$2$)] $e_i S\cd v_\chi^+ = 0  \fa i\in I$ and $\chi\in H_{\gamma,D}$; 
\item[(S$3$)] The highest degree PBW monomial of $S$ is a product of $f_i^{Dg_i}\fa i \in I$.
\endenu
\enddefn
\lem\label{lem:S exists}
Let $\gamma\in\Phi_{\text{re}}^+$ and $D\in\ZZ_{>0}$.
There is a Zariski-dense subset $\Tt$ of the hyperplane $H_{\gamma,D}$ such that for each $\chi\in\Tt$ there exists an element $S^\chi=S^\chi(\gamma,D)$ in $U(\fn_\CC^-)$ satisfying (S$1'$)-(S$3'$) below:
\enu
\item[(S$1'$)] $S^\chi = \sum\limits_{\pi\in\mathcal{P}} F_\pi Q_\pi^\chi$   where wt$(\pi)=D\gamma$  and $Q_\pi^\chi\in\CC$;
\item[(S$2'$)] $S^\chi\cd v^+_\chi$ is a singular vector in $M(\chi)_\CC$; 
\item[(S$3'$)] The highest degree PBW monomial of $S^\chi$ is a product of $f_i^{Dg_i}$ for all $i \in I$.
\endenu
Moreover, $\Tt$ satisfies a polynomial condition below:
\enu
\item[(S$4'$)] There are polynomials $Q_\pi \in\CC[h_0,\ldots, h_r]$ such that $Q_\pi(\chi) = Q_\pi^\chi$ for $\chi\in \Tt $.
\endenu
\endlem
\proof
The proof is exactly the same as in the finite case (see \cite[Sect. 4.13]{Hum08}). 
\endproof
\cor\label{cor:S exists}
For fixed $\gamma\in \Phi_{\text{re}}^+$ and $D\in \ZZ_{>0}$, there exists a (generic) Shapovalov element $S=S(\gamma,D)$.
Moreover, $S$ is unique  modulo the left $U(\fb_\CC^-)$-ideal $J(\gamma,D)= U(\fb^-_\CC)(h_\gamma +\rho(h_\gamma)-D)$.
\endcor
\proof
Let $S = \sum_\pi F_\pi Q_\pi $ where the polynomials $Q_\pi$ are defined in Lemma \ref{lem:S exists}. For each $\chi\in\Tt$ the evaluation $S(\chi)$ equals to the element $S^\chi$ constructed in Lemma \ref{lem:S exists}. 
By the density of $\Tt$, the evaluation $S(\chi)$ also satisfies Conditions (S$1'$) - (S$3'$) for each $\chi\in H_{\gamma,D}$ and hence $S$ satisfies Conditions (S$1$) - (S$3$).
\endproof
\subsection{Integral Shapovalov Elements}\label{sec:Int Shap elt}
The above-mentioned generic Shapovalov elements give rise to nonzero homomorphisms between Verma modules in characteristic 0. 
Franklin (\cite{Fr81}) introduced the integral Shapovalov elements for finite types and showed that such elements give rise to nonzero homomorphisms between the $\ZZ$-forms of Verma modules for finite types. 

We generalize his idea to the affine types and showed that the integral Shapovalov elements give rise to nonzero homomorphisms between the $\ZZ$-forms of Verma modules for affine types. In this section we mainly follow the outline given in \cite{Fr88}, while some adaptations are needed for the affine type.
\defn
An element $Z\in U(\fb_\CC^-)$ is called an \textit{integral Shapovalov element} 
(with respect to a positive real root $\gamma = \sum_{i\in I} g_i \af_i, g_i \in \ZZ$; and a positive integer $D$) 
if the conditions (Z1)--(Z3) below hold:
\enu
\item[(Z1)] $Z = \sum\limits_{\pi\in\mathcal{P}} F_\pi Q_\pi\text{ where wt}(\pi)=D\gamma\text{ and }Q_\pi\in\ZZ[h_0, \ldots, h_r]$;
\item[(Z2)] $
e_i^{(n)}Z \in I(\gamma,D)_\ZZ = \sum\limits_{\substack{j\in I\\ m\geq 1}}U_\ZZ e_j^{(m)} + 
U(\fb_\ZZ^-) (h_\gamma+\rho(h_\gamma) -D)
\text{ for } i\in I,n>0$;
\item[(Z3)] The highest degree PBW monomial in $Z$ is a product of $f_i^{Dg_i} \fa i \in I$.
\endenu
\enddefn
We will see that  for each $\gamma \in \Phi_{\text{re}}^+$ and $D\in \ZZ_{> 0}$, the integral Shapovalov element does exist and is unique (cf. \ref{cor:Z exists}), which we call  $Z(\gamma, D)$.  Here we prove the existence of $Z(\gamma,D)$ and its ``evaluation'' $Z^\chi(\gamma,D)$ simultaneously via an induction on the height of $\gamma$ as follows:
\enu
\item The initial case follows from an explicit construction.
\item Assuming the existence of $Z(\beta,D)$ for some $\beta = s_\af(\gamma) <\gamma$, 
we construct explicitly a dense subset $\Tt_1 = \Tt_1(\af)$ of $H_{\gamma,D}$ (cf. \eqref{eq:Tt1} ) such that 
for any $\chi \in \Tt_1$,  $Z^\chi$ exists.
\item The existence of  $Z^\chi$ in a dense subset $\Tt_1$ of $H_{\gamma,D}$ leads to the existence of $Z(\gamma, D)$.
\endenu
\thm \label{thm:gShap}
Let $\gamma = \sum_{i\in I} g_i \af_i \in\Phi_{\text{re}}^+$ and $D\in\ZZ_{>0}$. 
There is a dense subset $\Tt\subset H_{\gamma,D}$ such that for each $\chi\in \Tt$ there exists $Z^\chi=Z^\chi(\gamma,D) \in U(\fn_\CC^-)$ satisfying (Z$1'$)--(Z$3'$) below:
\enu
\item[(Z$1'$)] $Z^\chi = \sum\limits_{\pi\in\mathcal{P}} F_\pi Q_\pi^\chi\text{ where wt}(\pi)=D\gamma\text{ and }Q_\pi^\chi\in\ZZ$;
\item[(Z$2'$)] $
e_i^{(n)}Z^\chi \in I(\chi)_\ZZ$ (cf. \eqref{eq:IZ}) for $i \in I, n\in \ZZ_{>0};$  
\item[(Z$3'$)] The highest degree PBW monomial in $Z^\chi$ is a product of $f_i^{Dg_i}\fa i \in I$;
\item[(Z$4'$)] There are polynomials $Q_\pi \in\ZZ[h_0,\ldots, h_r]$ such that $Q_\pi(\chi) = Q_\pi^\chi \text{ for } \chi\in \Tt$.
\endenu 
Assuming Theorem \ref{thm:gShap}, one can prove the existence and uniqueness immediately.
\cor\label{cor:Z exists}
For each $\gamma\in\Phi_{\text{re}}^+$ and $D\in\ZZ_{>0}$ there exists an integral Shapovalov element 
$
Z=Z(\gamma,D) \in U(\fb_\ZZ^-).
$
Moreover, $Z$ is unique modulo the left $U(\fb_\ZZ^-)$-ideal 
\eq\label{eq:JZ}
J(\gamma,D)_\ZZ 
= U(\fb^-_\ZZ)(h_\gamma +\rho(h_\gamma)-D).
\endeq
\endcor
\proof
The existence part is similar to Corollary \ref{cor:S exists}.
The uniqueness part follows from that
$
e_i^{(n)}J(\gamma,D)_\ZZ \subset U(\fb^-_\ZZ) e_i^{(n)} +  U(\fb_\ZZ^-) (h_\gamma +\rho(h_\gamma)-D) \subset I(\gamma,D)_\ZZ.
$
\endproof

Before proceeding to the proof of Theorem \ref{thm:gShap} we need some lemmas. 

\lem\label{lem:gamma path}
If $\gamma\in\Phi_{\text{re}}^+$ then there are distinct positive roots $\{\gamma_i\}_{i=0}^n$ such that
\enu
\item $\gamma_0 = \gamma$.
\item $\gamma_i = s_{\epsilon_i}(\gamma_{i-1}) < \gamma_{i-1}$ for some simple root $\epsilon_i$.
\item $\gamma_n = \af_\eta$ for some simple root $\af_\eta$ occurring (cf. \S\ref{pre:root}) in $\gamma$.
\endenu
\endlem
\proof
This is done by induction on the height of $\gamma$.
\endproof
\rmk\label{rmk:diag}
For any $\ld\in H_{\gamma, Mp^e+D}$, we have a commutative diagram of weights via dot-actions:
\[
\bD{
\ld=\ld_0 & \rLine^{\epsilon_1} & \ld_1 &\rLine^{\epsilon_2} &\cds&\rLine^{\epsilon_{n-1}} &\ld_{n-1}\\
\uTo^{\gamma=\gamma_0}&&\uTo^{\gamma_1}&&\cds&&\uTo^{\gamma_{n-1}}\\
\mu=\mu_0 & \rLine_{\epsilon_1} & \mu_1 &\rLine_{\epsilon_2} &\cds&\rLine_{\epsilon_{n-1}} &\mu_{n-1}
}
\]
where $\ld_i= s_{\epsilon_i}\cd\ld_{i-1}, \mu_i= s_{\epsilon_i}\cd\mu_{i-1}=\ld_i-D\gamma$, and $\mu=\ld-D\gamma$.

Note that we reserve the freedom to choose a suitable dense subset of $H_{\gamma, Mp^e+D}$ to determine the order of the $\ld_i$'s and $\mu_i$'s, and hence the direction of the corresponding Verma module inclusions.
\endrmk
In the rest of this section we fix  $\gamma \in \Phi_{\text{re}}^+, D\in\ZZ_{>0}$,  $\chi \in H_{\gamma,D} \cap X$, $\nu = s_\af \cdot \chi$ and $\beta = s_\af(\gamma) < \gamma$ for  a simple root $\af$. We have 
\[
\ba{{l}
\nu = \chi + q\af ~\textup{ where }~q =-\g{\chi+\rho,\af^\vee}\in\ZZ, \textup{ and}\\
\beta = \gamma - b\af ~\textup{ where }~ b =  \g{\gamma, \af^\vee} \in \ZZ_{>0}.
}
\]
There are three possible commutative diagrams of weights via dot-actions:
\[
\ba{{ccc}
\textup{Case }1: q>0
&
\textup{Case }2:  0>q>-Db
&
\textup{Case }3:  q < - Db.
\\
\bD{
\chi&\rTo^{s_\af}&\nu\\
\uTo^{s_\gamma}&&\uTo_{s_\beta} \\
\chi-D\gamma&\rTo_{s_\af}&\nu-D\beta
}
\quad
&
\bD{
\chi&\lTo^{s_\af}&\nu\\
\uTo^{s_\gamma}&&\uTo_{s_\beta} \\
\chi-D\gamma&\rTo_{s_\af}&\nu-D\beta
}
\quad
&
\bD{
\chi&\lTo^{s_\af}&\nu\\
\uTo^{s_\gamma}&&\uTo_{s_\beta} \\
\chi-D\gamma&\lTo_{s_\af}&\nu-D\beta
}
}
\]
The first case will be used in proving Theorem \ref{thm:gShap}; the second case shows up in the proof of Lemma \ref{lem:CuZ}; although one can formulate a similar statement for the third case, it does not play a role in this article. Let
\eq\label{eq:Tt1}
\Tt_1 = \Tt_1(\gamma,D,\af) = H_{\gamma,D}\cap X\cap H_{\af,0}^-,
\endeq 
where $H_{\af,m}^-=\{\chi \in\fh^*_\RR\mid  \g{\chi+\rho,\af^\vee}<m\}$.
In other words, we have $q, q+Db \in \ZZ_{>0}$ for each $\chi \in \Tt_1$. 
We further let
\eq\label{eq:Tt2}
\Tt_2 = \Tt_2(\gamma,D,\af) = H_{\gamma,D}\cap X\cap H_{\af,0}^+ \cap H_{\af,-Db}^-,
\endeq
where $H_{\af,m}^+=\{\chi \in\fh^*_\RR\mid  \g{\chi+\rho,\af^\vee}>m\} \fa m\in\ZZ$.
In contrast, we have $q\in \ZZ_{<0}$ and $q+Db \in \ZZ_{>0}$ for each $\chi \in \Tt_2$.
\lem\label{lem:square}
Assume that $Z=Z(\beta,D)$ exists, then each $\chi\in \Tt_1$ (cf. \eqref{eq:Tt1}) has a unique (depending only on the choice of $Z$) element $Z^\chi = Z^\chi(\gamma,D,\af) \in U(\fn_\ZZ^-)_{-D\gamma}$ making the following diagram of Verma module homomorphisms commute.
\[
\bD{
M(\chi)_\CC&\rTo^{f_\af^q}&M(\nu)_\CC\\
\uTo^{Z^\chi}&&\uTo_{Z(\nu)} \\
M(\chi-D\gamma)_\CC&\rTo_{f_\af^{q+Db}}&M(\nu-D\beta)_\CC
}
\]
Here $Z(\nu)$ is the evaluation of $Z$ at $\nu$ (cf. \S\ref{pre:PBW}).
In particular, the element $Z^\chi$ is determined by the Shapovalov equation
\[
 Z^\chi f_\af^q= f_\af^{q+Db} Z(\nu) .
\]
\endlem
\proof
A direct calculation shows that $\nu-D\beta = s_\beta\cd\nu = s_\af \cd(\chi-D\gamma)$, and hence by the strong linkage principle (Proposition \ref{prop:lp})  we have a commuting diagram of nonzero homomorphisms below.
\[
\bD{
M(\chi)_\CC&\rTo^{i_3}&M(\nu)_\CC\\
\uTo^{i_1}&&\uTo_{i_4} \\
M(\chi-D\gamma)_\CC&\rTo_{i_2}&M(\nu-D\beta)_\CC
}
\]

Here $i_k$ is the left multiplication of a nonzero element $u_k \in U(\fn_\CC^-)$ for $k = 1,\ldots, 4.$ Since $f_\af^q \in \CC u_3, f_\af^{q+Db} \in \CC u_2$ and by assumption that $Z(\nu) \in \CC u_4$, so there is a unique element $Z^\chi \in \CC u_1$ such that the Shapovalov equation holds.
\endproof
Similarly, one can prove the counterpart of Lemma \ref{lem:square}.
\begin{lem}\label{lem:square b}
Assume that $Z=Z(\beta,D)$ exists, then each $\chi\in \Tt_2$ (cf. \eqref{eq:Tt2}) has a unique (depending only on the choice of $Z$) element  $Z^\chi = Z^\chi(\gamma,D,\af) \in U(\fn_\ZZ^-)_{-D\gamma}$ such that the following diagram of Verma module homomorphisms commutes.
\[
\bD{
M(\chi)_\CC&\lTo^{f_\af^{-q}}&M(\nu)_\CC\\
\uTo^{Z^\chi}&&\uTo_{Z(\nu)} \\
M(\chi-D\gamma)_\CC&\rTo_{f_\af^{q+Db}}&M(\nu-D\beta)_\CC
}
\]
In particular, the element $Z^\chi$ is determined by the Shapovalov equation
\[
 Z^\chi =f_\af^{q+Db} Z(\nu)  f_\af^{-q}.
\]
\end{lem}
Now we are finally in a position to prove Theorem \ref{thm:gShap}.
\proof[Proof of Theorem \ref{thm:gShap}]
The proof is done by induction on the height of $\gamma$. 
For the initial case $\gamma=\af_\eta$, we set $\Tt= H_{\gamma,D}\cap X$ and $Z^\chi = f_\eta^D$ so that the conditions (Z$1'$), (Z$3'$) and (Z$4'$) follow immediately from definition. Finally (Z$2'$) follows from the useful formula of Kostant below.
\[
e_i^{(n)}f_\eta^{(D)} = \delta_{i,\eta} \sum\limits_{k=0}^{\min(n,D)}f_i^{(D-k)} \ts{h_i-n-D+2k \choose k }e_i^{(n-k)}.
\]

For the inductive case ht$\gamma>1$, we apply Lemma \ref{lem:square} with $\af = \epsilon_1, \Tt = \Tt_1$
and obtain a Shapovalov equation
$
 Z^\chi f_\af^q= f_\af^{q+Db} Z(\nu) .
$
Let $Z(\nu) = \sum\limits_\omega F_\omega Q'_\omega$ where $Q'_\omega\in\ZZ$. 
Solving this leads to an explicit expression:
\[
\bA{
Z^\chi f_\af^q
&= \sum_{\omega\in\mathcal{P}} f_\af^{q+Db} F_\omega Q'_\omega
\\
&= \sum_{\substack{\pi \in \mathcal{P}\\  \pi(\af) = 0}} 
\sum_{i\in\ZZ} F_\pi f_\af^{q+Db-i} \left(\sum_{\omega\in\mathcal{P}}  C_{\omega,\pi,i}(q+Db) Q'_\omega\right).
\\
}
\]
Therefore, $Z^\chi = \sum_\pi F_\pi Q_\pi^\chi$ 
where
\[
Q_\pi^\chi =\sum_{\omega\in\mathcal{P}} C_\omega Q'_\omega \textup{ for some }C_\omega \in \ZZ.
\]
It is routine to check that the element $Z^\chi$ obtained this way satisfies (Z$1'$)--(Z$3'$).
%
%
%
%
For (Z$4'$), by the inductive hypothesis for each $\pi$  there is a polynomial $P_\pi \in \ZZ[h_0, \ldots, h_r]$ such that
$Q'_\pi = P_\pi(\nu) 
.$ Let $\af =\af_k$. Note that for each $i\in I$ we have
\[
h_i(\nu) = \g{\nu,\af_i^\vee} = \g{\chi,\af_i^\vee} - \g{\chi+\rho,\af^\vee}\g{\af,\af_i^\vee} = \bp{h_i-a_{ik}(h_k+1)}(\chi).
\]
So there is a polynomial $P'_\pi \in \ZZ[h_0,\ldots,h_r]$, obtained from $P_\pi$ by replacing $h_i$ by $h_i-a_{ik}(h_k+1)$ for each $i$, such that
$P'_\pi(\chi)  = P_\pi(\nu) = Q'_\pi$
 and hence by defining $Q_\pi = \sum_{\omega\in\mathcal{P}} C_\omega P'_\pi$  one has
\[
Q_\pi(\chi) = \sum_{\omega\in\mathcal{P}} C_\omega P'_\pi(\chi) 
= \sum_{\omega\in\mathcal{P}} C_\omega Q'_\pi = Q_\pi^\chi.
\]
\endproof

\cor
Assume that $\ld \in H_{\gamma,D}$ and $\mu = \ld-D\gamma$.
There is a nonzero homomorphism $M(\mu)_\ZZ \rw M(\ld)_\ZZ$ given by
\[
v_\mu^+ \mapsto Z v_\ld^+,
\] 
where $Z = Z(\gamma,D)$ is the integral Shapovalov element.
\endcor
\exa \label{exa}
Consider the affine Kac-Moody algebra of type $\widetilde{A}_1$ with
fixed positive real root $\gamma = \af_0+\delta$ and $D = 1$.
Then the sequence of positive roots described in Lemma \ref{lem:gamma path} is given by 
\[
\af_0+\delta = \gamma_0 > \gamma_1 = s_{\af_0}(\gamma_0) = \af_1.
\]
The dense subset $\Tt\subset H_{\gamma,D}$  in Theorem \ref{thm:gShap} is given by 
\[
\bA{
\Tt
&= \left\{\chi =\xi_0\varpi_0-2(\xi_0+1)\varpi_1+\xi\delta \mid  \xi_0\in \ZZ_{<0} \right\}\\
&= \left\{-(q+1)\varpi_0+2q\varpi_1+\xi\delta \mid  q\in \ZZ_{\geq0} \right\},
}
\] 
where $q=-\g{\chi+\rho,\af_0^\vee}$.
For each weight $\chi\in \Tt$, we have a Shapovalov equation
\[
Z^\chi f_0^q = f_0^{q+2} f_1.
\] 
Solving this, we get $Z^\chi = f_1f_0^2 + f_{\delta,1} f_0 (q+2) + f_\gamma (q+2)(q+1)$ as in Example \ref{exaC}. 
Note that evaluation at $\chi = -(q+1)\varpi_0+2q\varpi_1+\xi\delta$ replaces $h_0$ by $-(q+1)$, equivalently, $q+2 = -h_0(\chi)+1$.
Hence
it leads to  an integral Shapovalov element
\[
Z = f_1f_0^2 - f_{\delta,1} f_0 (h_0-1) + f_\gamma (h_0^2 - h_0).
\]
\endexa
%

\section{Homomorphisms between Verma modules in characteristic $p$}  \label{sec:Verma}
For finite types, Franklin (\cite{Fr88}) showed that if the corresponding positive root for an integral Shapovalov element is not ``exceptional'' , 
there is a ``modified'' Shapovalov element that gives rise to a nonzero homomorphism between Verma modules in characteristic $p$, without any restriction. 
As for the five ``exceptional'' roots, additional restrictions on the distance between the two highest weights are required. 

In  this section we generalize Franklin's idea to affine types.
In our viewpoint, the restriction for ``exceptional'' roots  is in fact a generic phenomenon.
We first show that under a restriction on the distance, integral Shapovalov elements give rise to nonzero homomorphisms between Verma modules in characteristic $p$. 
We then show that the restriction is the coarsest if the corresponding positive real root is ``good''.
 
\subsection{$\eta$-goodness}
\defn\label{defn:eta-good}
For a fixed prime number $p$ and $\eta\in I$, a positive real root $\gamma \in \Phi_{re}^+$ is called \textit{$\eta$-good} if 
the multiplicity of $h_\eta$ occurring in $h_\gamma$ (cf. \eqref{eq:hr}) is coprime to $p$, and 
the simple root $\af_\eta$ has the same length as $\gamma$. 

In other words, assume that $h_\gamma = \sum g^\vee_i h_i$ for some $g_i^\vee \in \ZZ$. $\gamma$ is $\eta$-good if gcd$(g^\vee_\eta,p )= 1$ and $(\af_\eta,\af_\eta) = (\gamma,\gamma)$.
\enddefn
\lem\label{lem:good}
If $\gamma$ is $\eta$-good then
\[
\ZZ[h_0,\ldots,h_r]\otimes\ZZ_{(p)} \simeq \ZZ[h_0,\ldots,\widehat{h}_\eta,\ldots,h_r, h_\gamma + \rho(h_\gamma) -D]\otimes\ZZ_{(p)}.
\]
Here $\widehat{h}_\eta$ represents the omission of the variable and $\ZZ_{(p)}=\{r/q \in \QQ \mid  \gcd(p,q) = 1\}$ is the localization of $\ZZ$ at $p\ZZ$.
\endlem
\proof~
It follows from that $h_\eta = \frac{1}{g_\eta^\vee}(h_\gamma - \sum_{i\neq \eta} g_i^\vee h_i)$ and $p\not|  ~g^\vee_\eta$.
\endproof
Now we define an ``$h_\eta$-avoiding'' Shapovalov-type element by modifying $Z$ as follows.
\lem\label{lem:repn}
Let $Z = Z(\gamma,D) = \sum_\pi F_\pi Q_\pi$ with $Q_\pi \in \ZZ[h_0, \ldots, h_r]$ for some $D\in\ZZ_{>0}$ and  $\gamma\in\Phi_{\text{re}}^+$  
 such that $\af_\eta$ is occurring in $\gamma$,
and let $N = \max \{\deg_{h_\eta} Q_{\pi} ~|~ \pi\in \mathcal{P}\}$.
There is an element $Z_\eta = Z_\eta(\gamma,D) \in U(\fb^-_\ZZ)$ satisfying the following conditions (G1)-(G3):
\enu
\item[(G1)] $Z_\eta = \sum\limits_{\pi\in\mathcal{P}} F_\pi Q_{\eta,\pi}\text{ where wt}(\pi)=D\gamma\text{ and }Q_{\eta,\pi}\in\ZZ[ h_0, \ldots, \widehat{h}_\eta, \ldots, h_r]$;
\item[(G2)] $e_i^{(n)}Z_\eta \in I(\gamma,D)_\ZZ  \fa i\in I$ and $n>0$;
\item[(G3)] The highest degree PBW monomial in $Z_\eta$ is a product of $f_i^{Dg_i} \fa i \in I$ times $(g^\vee_\eta)^N$.
\endenu
\endlem
\proof
By definition, 
\[
(g^\vee_\eta)^N Z = \sum\limits_{\pi\in\mathcal{P}} F_\pi  (g^\vee_\eta)^N Q_\pi = \sum\limits_{\pi\in\mathcal{P}} F_\pi Q'_\pi,
\]
where $Q'_\pi \in \ZZ[h_0, \ldots, \widehat{h}_\eta, \ldots, h_r, h_{\gamma}+\rho(h_\gamma)-D]$. Equivalently,
\[
(g^\vee_\eta)^N Z = \sum_{m=0}^N u_{\eta,m} (h_\gamma+\rho(h_\gamma) -D)^m,
\]
for some $u_{\eta,m} \in U(\fb^-_\ZZ)$. 
Finally, we define 
\[
Z_\eta = (g^\vee_\eta)^N Z -\sum_{m = 1}^N u_{\eta,m} (h_\gamma + \rho(h_\gamma) - D)^m
=\sum\limits_{\pi\in\mathcal{P}} F_\pi Q_{\eta,\pi},
\] 
where $Q_{\eta,\pi} \in \ZZ[h_0, \ldots, \widehat{h}_\eta, \ldots, h_r]$ as required in (G1).
By construction (G2) is satisfied. (G3) follows from the fact that 
$e_i^{(n)} (h_\gamma + \rho(h_\gamma) -D) \in U_\ZZ e_i^{(n)} \subset I(\gamma,D)_\ZZ$ 
and hence
\[
e_i^{(n)}Z_\eta = (g^\vee_\eta)^N e_i^{(n)} Z 
-\sum_{m = 1}^N 
u_{\eta,m} e_i^{(n)} (h_\gamma + \rho(h_\gamma) - D)^m
 \in  I(\gamma,D)_\ZZ.
\]
\endproof
\exa\label{exa4}
Recall that in Example \ref{exa} we have $\gamma = 2\af_0 + \af_1$ and hence $h_\gamma = 2h_0 + h_1$.
Therefore $\gamma$ is $0$-good for any odd prime and is $1$-good for any prime.
We have already an $h_1$-avoiding Shapovalov element
\[
Z_1  = Z = f_1f_0^2 - f_{\delta,1} f_0 (h_0-1) + f_\gamma (h_0^2 - h_0).
\] 
As for $\eta = 0$, to achieve an $h_0$-avoiding Shapovalov element we first replace $h_0$ by $\frac{1}{2}(h_\gamma - h_1)$ and then make a suitable scaling. Thus
\[
g_0^N Z =  2^2Z =4f_1f_0^2 - 2f_{\delta,1} f_0 (h_\gamma - h_1 - 2) + f_\gamma (h_\gamma - h_1)(h_\gamma - h_1-2).
\]
Note that $h_\gamma + \rho(h_\gamma) - D = h_\gamma +2$.  Setting $h_\gamma = -2$, we obtain
\[
Z_0 =  4f_1f_0^2 + 2f_{\delta,1} f_0 ( h_1 + 4) + f_\gamma (h_1+2)(h_1+4).
\]
\endexa
\subsection{Projection onto $U(\fb_\CC^-)$}
Let $P_{\fb^-}:U\rw U(\fb_\CC^-)$ be the projection corresponding to the direct sum decomposition $U = \sum_{\af\in\Phi^+}U e_\af \oplus U(\fb_\CC^-)$. 
It is easy to see that for each $\ld\in\fh_\RR^*, u \in U$ we have
\[
uv_\ld^+ = P_{\fb^-}(u) v_\ld^+.
\] 
In this section we study $e_i^{(n)} Z_\eta v_\ld^+$ for the $h_\eta$-avoiding integral Shapovalov element $Z_\eta$ we constructed in Lemma \ref{lem:repn}.
\lem\label{lem:proj}~
\enu
\item[(1)]
Assume that $n\in\ZZ_{>0}, i \in I$. Let $\Psi_{\pi,\omega} =\Psi_{\pi,\omega}(i,n)\in U(\fh_\CC)$ be such that
\[
P_{\fb^-}\left(e_i^{(n)}F_\pi\right) 
= \sum_{\omega\in\mathcal{P}} F_\omega \Psi_{\pi,\omega}
.
\]
Then $\Psi_{\pi,\omega} \in \ZZ[h_i]$ and $\deg_{h_i}\Psi_{\pi,\omega} \leq \pi(\af_i)$.
\item[(2)]
Let $D\in\ZZ_{>0}, \gamma= \sum_{i\in I} g_i\af_i \in \Phi_{\text{re}}^+$ with $g_\eta \neq 0$,
$Z_\eta=Z_\eta(\gamma,D) = \sum\limits_{\pi\in\mathcal{P}} F_\pi Q_{\eta,\pi}$ where $Q_{\eta,\pi} \in \ZZ[h_0, \ldots,\widehat{h}_\eta,\ldots, h_r]$ as in Lemma \ref{lem:repn},
and let $\Psi'_{\omega} \in U(\fh_\CC)$ be such that 
\[
P_{\fb^-}\left(e_i^{(n)}Z_\eta \right) 
= \sum_{\omega\in\mathcal{P}} F_\omega \Psi'_{\omega}.
\]
Then  $\Psi'_{\omega} \in \ZZ[h_0,\ldots,h_r]$ and $\deg_{h_\eta}\Psi'_{\omega} \leq Dg_\eta$.
\endenu
\endlem
\proof
The first part is exactly the same as in the finite case \cite[Lemma 4.3]{Fr88}.
Also,
\[
\bA{
\deg_{h_\eta} \Psi'_\omega 
&= \max\left\{\deg_{h_\eta}\Psi_{\pi,\omega} + \deg_{h_\eta} Q_{\eta,\pi} ~\middle|~  \textup{wt}(\pi) = D\gamma \right\}
\\
&\leq \delta_{i,\eta}\max\left\{\pi(\af_\eta) \mid  \textup{wt}(\pi) = D\gamma \right\} 
\leq Dg_\eta. 
}
\]
This proves (2).
\endproof
Combining Lemma \ref{lem:good} and Lemma \ref{lem:proj}, we can express $P_{\fb^-}\left(e_i^{(n)}Z_\eta\right) $ as
\eq\label{eq:PbeZ}
P_{\fb^-}\left(e_i^{(n)}Z_\eta\right)  
= \sum_{m=0}^{Dg_\eta} u_m \ts{h_\gamma + \rho(h_\gamma) -D \choose m},
\endeq
for some $u_m \in U(\fn_\ZZ^-) \otimes \ZZ[h_0,\ldots,\widehat{h}_\eta,\ldots, h_r]\otimes \ZZ_{(p)}$.
\prop\label{prop:PbeS}
Let $Z_\eta=Z_\eta(\gamma,D)$. 
If  $\gamma=\sum g_i\af_i \in \Phi_{\text{re}}^+$ is $\eta$-good then 
$h_\gamma +\rho(h_\gamma) - D$ divides $P_{\fb^-}\left(e_i^{(n)}Z_\eta\right)$, that is, $u_0$ defined in \eqref{eq:PbeZ} equals to 0. 
\endprop
\proof
For each $\chi\in H_{\gamma,D}$ we have  $e_i^{(n)}Z_\eta v^+_\chi = 0$ and hence $ P_{\fb^-}\left(e_i^{(n)}Z_\eta\right)v_\chi^+ =0$.
Thus $h_\gamma +\rho(\gamma) - D$ divides $P_{\fb^-}\left(e_i^{(n)}Z_\eta \right)$.
\endproof
\subsection{Homomorphisms between Verma modules}
Now we are in a position to construct the homomorphisms in detail.
Note that the naive map $M(\mu)_\KK \rw M(\ld)_\KK$ sending $v_\mu^+\otimes 1_\KK$ to $Z_\eta v_\ld^+\otimes 1_\KK$ is not necessarily nonzero since $p$ may divide $Z_\eta(\ld)$ in $U_\ZZ$. 
In our construction we send $v_\mu^+\otimes 1_\KK$ to  $(Z_\eta v_\ld^+)/p^f\otimes 1_\KK$ instead, where $p^f$ is the highest $p$-power dividing $Z_\eta(\ld)$ in $U_\ZZ$.
\prop\label{thm:PbeS}
If $N\in\ZZ_{>0}$ divides $Z_\eta(\ld)$ in $U_\ZZ$ then $N$ also divides each $u_m$ defined in \eqref{eq:PbeZ}.
\endprop
\proof
%
 It is exactly the same as in the finite case \cite[Proposition 5.1]{Fr88}.
\endproof
\cor
Assume that $\mu\uw_\gamma \ld$ (cf. \S\ref{pre:low}) and $\gamma=\sum g_i \af_i \in \Phi_{\text{re}}^+$ is $\eta$-good.
If $Dg_\eta < p^e$ then there is a nonzero homomorphism $M(\mu)_\KK \rw M(\ld)_\KK$.
\endcor
\proof
Let $p^f$ be the highest $p$-power dividing $Z_\eta(\ld)$ in $U_\ZZ$.
Applying Lemma \ref{prop:PbeS} and Proposition \ref{thm:PbeS} with $N = p^f$, we have
\[
e_i^{(n)} \frac{Z_\eta(\ld)}{p^f} v_\ld^+
= \sum_{m=1}^{Dg_\eta} \frac{u_m}{p^f} 
\ts{\g{\ld+\rho,\gamma^\vee} -D \choose m} v_\ld^+.
\]
Since $Dg_\eta < p^e$, each $m <p^e$ and hence each ${\g{\ld+\rho,\gamma^\vee}-D \choose m} = {Mp^e \choose m} \in p\ZZ$. 
Therefore $e_i^{(n)}  \frac{Z_\eta(\ld)}{p^f} v_\ld^+\otimes 1_\KK = \overline{0}$ and the map $M(\mu)_\KK \rw M(\ld)_\KK$ sending
$v_\mu^+\otimes 1_\KK$ to $\frac{Z_\eta(\ld)}{p^f} v_\ld^+\otimes 1_\KK$ is a homomorphism.
\endproof
\section{Homomorphisms between Weyl modules}  \label{sec:Weyl}
Now we start to deal with the Weyl modules using reduction modulo $p$.
To construct an analogous homomorphism as in Section \ref{sec:Verma}, we have to make sure that  $Z_\eta  \overline{v}_\ld^+$ is nonzero in $L(\ld)_\ZZ$, which is not obvious.
This is done by using the contravariant form and the Shapovalov factor formula (= Lemma \ref{lem:CuZ}).
Finally we prove the existence of nonzero homomorphisms between Weyl modules whose highest weights are $\gamma$-linked if the two highest weights are close enough.  In this section we mainly follow the outline given in \cite{Fr81}, while some adaptations are needed for the affine type.
%
\subsection{Contravariant Forms}
In this section we review some basic properties  of the contravariant form mentioned in \cite[Sect. 3.14, 3.15]{Hum08} and \cite[\S8]{Fr81}.
A symmetric bilinear form $C:U\times U \rw U(\fh_\CC)$ is called \textit{contravariant} if
\[
C(u\cdot v,v') = C(v,\tau(u)\cdot v') \fa u\in U, v, v'\in U.
\]
\prop\label{prop:contra}~ 
\enu
\item[(a)] Let $P_\fh:U \rw U(\fh_\CC)$ be the projection. Then $C(u,v)= P_\fh(\tau(u)v)$ is a nonzero contravariant form.
\item[(b)] 
The weight spaces are orthogonal to each other with respect to $C$, i.e., $C(U_\ld, U_\mu) = 0$ if $\ld\neq \mu$.
\item[(c)] 
Let $C^\ld:U\times U \rw \CC$ be the evaluation of $C$ at $\ld$. Then
\[
\Rad C^\ld = \{u\in U ~|~ u\overline{v}^+_\ld=0\}.
\]
\item[(d)] 
$C$ is non-degenerate on $U(\fb^-_\CC)$.
\item[(e)] 
$\deg C(F_\pi, F_\omega) \leq \min (\deg F_\pi, \deg F_\omega)$.
\endenu
\endprop
\subsection{Shapovalov Factor Formula}
Recall that for arbitrary Kac-Moody algebras, the Shapovalov determinant formula \cite[Theorem 1]{KK79} computes the determinant of the contravariant form on $U_{-D\gamma}$ as follows.
\[
\text{det} C\mid _{U_{-D\gamma}} = 
\prod_{\af\in \Phi_{\text{re}}^+} \prod_{n=1}^\infty \left(h_\af + \rho(h_\af) - n \right)^{P(D\gamma - n\af)} 
\prod_{\af\in \Phi_{\text{im}}^+} \left(h_\af + \rho(h_\af)  \right)^{\sum\limits_{n=1}^\infty P(D\gamma - n\af)}.
\]

Here $P$ is the Kostant partition function.
In this section we prove a variant that describes the common factor for the column corresponding to an integral Shapovalov element $Z(\gamma,D)$.
\lem \label{lem:CuZ}
Let $Z = Z(\gamma,D)$ for some $\gamma \in \Phi_{\text{re}}^+$ and $D\in \ZZ_{>0}$.
Define $\epsilon_i, \gamma_i$ as in Lemma \ref{lem:gamma path}.
Let $b_i \in\ZZ_{\geq0}$ be such that $\gamma = \sum b_i \epsilon_i$,
and let $\beta_i=s_{\epsilon_1}\cds s_{\epsilon_{i-1}}(\epsilon_i) \in \Phi_{\text{re}}^+$.
Then for any $u\in U$ we have
\[
C(u,Z) \in\CC \prod_{i=1}^n \prod_{j=1}^{Db_i} (h_{\beta_i} +\rho(h_{\beta_i}) - j) .
\]

\endlem
\proof
We will show that for fixed $i\in\{1,\ldots, n\}$ and $j\in\{1,\ldots,Db_i\}$ that $h_{\beta_i}+\rho(h_{\beta_i})-j$ divides $C(u,Z)$.
By Proposition \ref{prop:contra} one needs to check that
$
C^\ld(u,Z) = 0 \fa \ld \in H_{\beta_i,j}
$,
equivalently, $Z v_\ld^+ \in N(\ld)_\CC $.
So it remains to check that there is a dense subset $\Tt$ of $H_{\beta_i,j}$ such that
for any $\ld \in \Tt$, we have
\[
Z v_\ld^+ \in M(s_\af\cd\ld)_\CC, 
\quad \af \in \Phi_{\text{re}^+}.
\]
Now $q_i= - \g{\ld+\rho,\beta_i^\vee}=-j < 0$.
By assumption $q_i+Db_i =Db_i-j \geq 0$, and hence by Lemma \ref{lem:square b} there is a 
commutative diagram of Verma module inclusions
\eq\label{diag:A}
\bD{
M(\ld_{i-1})_\CC&\lTo^{f_{\epsilon_i}^{-q_i}}&M(\ld_i)_\CC \\
\uTo^{Z(\gamma_{i-1},D)(\ld_{i-1})}&&\uTo_{Z(\gamma_i,D)(\ld_i)}\\
M(\mu_{i-1})_\CC&\rTo^{f_{\epsilon_i}^{q_i+Db_i}}&M(\mu_i)_\CC 
}
\endeq

Recall that as in Remark \ref{rmk:diag} we can choose a dense subset of $H_{\gamma, MP^e+D}$ such that for any $\ld$ in the subset we have the following commutative diagram of weights.
\[
\bD{
\ld=\ld_0 & \lTo^{\epsilon_1} & \ld_1 &\lTo^{\epsilon_2} &\cds&\lTo^{\epsilon_{i-1}} &\ld_{i-1}\\
\uTo^{\gamma=\gamma_0}&&\uTo^{\gamma_1}&&&&\uTo^{\gamma_{i-1}}\\
\mu=\mu_0 & \lTo_{\epsilon_1} & \mu_1 &\lTo_{\epsilon_2} &\cds&\lTo_{\epsilon_{i-1}} &\mu_{i-1}
}
\]
By Lemma \ref{lem:square} the above diagram leads to a commutative diagram of Verma module inclusions.
\eq\label{diag:B}
\bD{
M(\ld_0)_\CC & \lTo^{f_{\epsilon_1}^{q_1}} & M(\ld_1)_\CC &\lTo^{f_{\epsilon_2}^{q_2}} &\cds&\lTo^{f_{\epsilon_{i-1}}^{q_{i-1}}} &M(\ld_{i-1})_\CC\\
\uTo_{Z(\gamma,D)}&&\uTo_{Z(\gamma_1,D)(\ld_1)}&&&&\uTo^{Z(\gamma_{i-1},D)(\ld_{i-1})}\\
M(\mu_0)_\CC & \lTo_{f_{\epsilon_1}^{q_1-Db_1}} & M(\mu_1)_\CC &\lTo_{f_{\epsilon_2}^{q_2-Db_2}} &\cds&\lTo_{f_{\epsilon_{i-1}}^{q_{i-1}-Db_{i-1}}} &M(\mu_{i-1})_\CC
}
\endeq
Finally, calculations show that if $i\neq n$ then for a suitable choice of dense subset of $H_{\gamma, MP^e+D}$ there is a commutative diagram of weights as the following with $\ld'_{i-1}=\ld_i$, $\ld'_k=s_{\epsilon_{k+1}}\cd \ld'_{k+1}$, $\mu'_{i-1}=\mu_i$, $\mu'_k=s_{\epsilon_{k+1}}\cd \mu'_{k+1}$, $\beta'_k = s_{\epsilon_{k+1}}\cds s_{\epsilon_{i-1}}(\epsilon_i)$, $\gamma'_{i-1}=\gamma_i$ and $\gamma'_k=s_{\epsilon_{k+1}}(\gamma_{k+1})$.
\[
\bD{
\ld_0 & \lTo^{\epsilon_1} & \ld_1 &\lTo^{\epsilon_2} &\cds&\lTo^{\epsilon_{i-1}} &\ld_{i-1}\\
\uTo_{\beta_i=\beta'_0}&&\uTo_{\beta'_1}&&&&\uTo^{\beta'_{i-1}=\epsilon_i}\\
\ld'_0 & \lTo^{\epsilon_1} & \ld'_1 &\lTo^{\epsilon_2} &\cds&\lTo^{\epsilon_{i-1}} &\ld'_{i-1}=\ld_i\\
\uTo_{\gamma_i=\gamma'_0}&&\uTo_{\gamma'_1}&&&&\uTo^{\gamma'_{i-1}=\gamma_i}\\
\mu'_0 & \lTo_{\epsilon_1} & \mu'_1 &\lTo_{\epsilon_2} &\cds&\lTo_{\epsilon_{i-1}} &\mu'_{i-1}=\mu_i\\
\uTo_{\beta_i=\beta'_0}&&\uTo_{\beta'_1}&&&&\uTo^{\beta'_{i-1}=\epsilon_i}\\
\mu_0 & \lTo_{\epsilon_1} & \mu_1 &\lTo_{\epsilon_2} &\cds&\lTo_{\epsilon_{i-1}} &\mu_{i-1}
}
\]
Now let $Q_k= \g{\gamma'_k,\epsilon_k^\vee}$. Again by Lemma \ref{lem:square}  the above diagram leads to a commutative diagram  of Verma module inclusions.
\eq\label{diag:C}
\bD{
M(\ld_0)_\CC & \lTo^{f_{\epsilon_1}^{q_1}} & M(\ld_1)_\CC &\lTo^{f_{\epsilon_2}^{q_2}} &\cds&\lTo^{f_{\epsilon_{i-1}}^{q_{i-1}}} &M(\ld_{i-1})_\CC\\
\uTo_{S(\beta_i,j)}&&\uTo&&&&\uTo^{f_{\epsilon_i}^j}\\
M(\ld'_0)_\CC & \lTo^{f_{\epsilon_1}^{q_1-jb_1}} & M(\ld'_1)_\CC &\lTo^{} &\cds&\lTo^{} &M(\ld'_{i-1})_\CC\\
\uTo_{Z(\gamma_i,D)(\ld'_0)}&&\uTo&&&&\uTo^{Z(\gamma'_{i-1},D)(\ld'_{i-1})}\\
M(\mu'_0)_\CC & \lTo^{f_{\epsilon_1}^{q_1-jb_1-DQ_1}} & M(\mu'_1)_\CC &\lTo_{} &\cds&\lTo_{} &M(\mu'_{i-1})_\CC\\
\uTo_{S(\beta_i,Db_i-j)(\mu'_0)}&&\uTo&&&&\uTo^{f_{\epsilon_i}^{Db_i-j}}\\
M(\mu_0)_\CC & \lTo^{f_{\epsilon_1}^{q_1-Db_1}} & M(\mu_1)_\CC &\lTo_{} &\cds&\lTo_{} &M(\mu_{i-1})_\CC
}
\endeq

Combining Diagrams \eqref{diag:A}, \eqref{diag:B} and \eqref{diag:C} 
and using the fact that $U$ has no zero divisors, 
we obtain another commutative diagram below:
\[
\bD{
M(\ld)_\CC&\lTo^{S(\beta_i,j)}&M(\ld'_0)_\CC \\
\uTo^{Z(\gamma,D)}&&\uTo_{Z(\gamma_i,D)(\ld'_0)}\\
M(\mu)_\CC&\lTo^{S(\beta_i,Db_i-j)(\mu'_0)}&M(\mu'_0)_\CC 
}
\]
This shows that $Z v^+_\ld = S(\beta_i,Db_i-j)(\mu'_0)Z(\gamma_i,D)(\ld'_0)S(\beta_i,j) v^+_\ld \in M(s_{\beta_i}\cd\ld)_\CC$.

For the special case $i=n$, Diagram  \eqref{diag:C}  is replaced by
\eq
\label{diag:C'}
\bD{
M(\ld_0)_\CC & \lTo^{f_{\epsilon_1}^{q_1}} & M(\ld_1)_\CC &\lTo^{f_{\epsilon_2}^{q_2}} &\cds&\lTo^{f_{\epsilon_{i-1}}^{q_{i-1}}} &M(\ld_{i-1})_\CC\\
\uTo_{S(\gamma,j)}&&\uTo&&&&\uTo^{f_{\epsilon_i}^j}\\
\ba{{r}M(\ld'_0)_\CC\\=M(\mu'_0)_\CC} & \lTo^{f_{\epsilon_1}^{q_1-jb_1}} & \ba{{r}M(\ld'_1)_\CC\\=M(\mu'_1)_\CC} &\lTo^{} &\cds&\lTo^{} &\ba{{r}M(\ld'_{i-1})_\CC\\=M(\mu'_{i-1})_\CC}\\
\uTo_{S(\gamma,D-j)(\mu'_0)}&&\uTo&&&&\uTo^{f_{\epsilon_i}^{Db_i-j}}\\
M(\mu_0)_\CC & \lTo^{f_{\epsilon_1}^{q_1-Db_1}} & M(\mu_1)_\CC &\lTo_{} &\cds&\lTo_{} &M(\mu_{i-1})_\CC
}.
\endeq
Similarly we have
$Z v^+_\ld \in M(s_{\beta_n}\cd\ld)_\CC = M(s_\gamma\cd\ld)_\CC$.

Note that the $\beta_i$'s are distinct and none of them are multiples of others, so the factors $h_{\beta_i} + \rho(h_{\beta_i}) - j$ are relatively prime.
Using Proposition \ref{prop:contra}(e), one can show that those are all the factors, and hence each $C(u,Z)$ is a scalar multiple of the product of these factors. 

\endproof
%
\cor
If $\gamma$ is $\eta$-good then $Z_\eta \overline{v}^+_\ld \neq \overline{0}$ in $L(\ld)_\ZZ$.
\endcor
\proof
Suppose that $Z_\eta \overline{v}^+_\ld = \overline{0}$ in $L(\ld)_\ZZ$, then $Z \overline{v}^+_\ld = \overline{0}$ in $L(\ld)_\ZZ$ as well. Therefore $C^\ld(u, Z(\ld)) = 0 \fa u \in U$ and hence  $C^\ld(u, Z) = 0 \fa u \in U$.
By Proposition \ref{prop:contra} we know that $C$ is non-degenerate on $U(\fb^-_\CC)$ and hence the fact $Z \neq 0$ implies $C(u,Z)\neq 0$ for some $u \in U$.

Applying Lemma \ref{lem:CuZ}, we have $\prod_{i=1}^n \prod_{j=1}^{Db_i} (h_{\beta_i} +\rho(h_{\beta_i}) - j)(\ld) = 0$, 
and hence 
$\g{\ld+\rho,\beta_i^\vee}=j$ for some $1\leq  i\leq n $ and $1 \leq j  \leq Db_i$.
Therefore both $j$ and $Db_i - j $ are non-negative so that $\ld_i$ and $\mu_i$ are on the opposite sides of the hyperplane $H_{\epsilon_i,0}$. 
This shows that they are in different chambers,  which is a contradiction.
\endproof 
\subsection{BGG resolution}
Kumar has generalized the strong BGG resolution over any Kac-Moody algebras. 
\prop \label{prop:BGG}
For each $\ld\in X^+$ there is an exact sequence of $\fg_\CC$-modules:
\[
\cds \rw \mc{C}_L \rw \cds \rw \mc{C}_1 \rw M(\ld)_\CC \rw L(\ld)_\CC \rw 0,
\]
where $\mc{C}_L= \mathop{\bigoplus\limits_{\substack{ w\in W\\ l(w) = L}}} M(w\cd\ld)_\CC$.
\endprop
\proof
See \cite[Theorem 3.20]{ku90}.
\endproof
Taking advantage of the BGG resolution, we describe bases of certain weight spaces of Verma modules and Weyl modules in characteristic 0.
\lem\label{lem:wt}
Assume that $\ld\in X^+, \beta=\sum m_i \af_i \in Q^+$ and $w\in W$.
If there is an $\eta \in I$ such that  $s_{\af_\eta}$ occurs in $w$ and $\g{\ld+\rho,\af_\eta^\vee}> m_\eta $, 
then
$\ld-\beta$ is not a weight of $M(w\cd\ld)_\CC$. 
In particular,
each $\mathcal{C}_L$ (cf. Proposition \ref{prop:BGG}) has weight space 
\[
(\mathcal{C}_L)_{\ld-\beta} =\bigoplus_{\substack{ w\in W_{I\backslash\{\eta\}}\\l(w)=L } } M(w\cd\ld)_{\CC,\ld-\beta},
\]
where  $W_{I\backslash\{\eta\}}$ is the parabolic subgroup of $W$ generated by $\{s_{\af_i}\}_{i \in I\backslash\{ \eta\}}$.
\endlem
\proof
Fix a reduced expression $w= s_{\af_{i_N}}\cds s_{\af_{i_1}}$  and let $\beta'_j = s_{\af_{i_N}} \cds s_{\af_{i_{j+1}}}(\af_{i_j})$.
Suppose that $\ld-\beta$ is a weight of $M(w\cd\ld)_\CC$, then
\[
\ld-\beta < w\cd\ld = \ld - \sum_{j=1}^N \g{\ld+\rho,\af_{i_j}^\vee} \beta'_j.
\]
Since $s_{\af_\eta}$ occurs in $w$, there is an integer $L$ which is the largest one satisfying $\af_\eta = \af_{i_L}$,  that is,
\[
\beta'_L = s_{\af_{i_N}} \cds s_{\af_{i_{L+1}}}(\af_\eta) =\af_\eta + \sum_{k\neq \eta} n_k \af_k \fsm n_k\geq 0.
\]
Hence,
\[
\ld - \beta < \ld - \g{\ld+\rho,\af_{i_L}^\vee} \beta'_L < \ld - \g{\ld+\rho,\af_\eta^\vee}\af_\eta,
\]
and therefore $\g{\ld+\rho,\af_\eta^\vee}< m_\eta$, which is a contradiction.
\endproof
Assume that $\ld\in X^+, w\in W$ with a fixed reduced expression $w = s_{\af_{i_N}}\cds s_{\af_{i_1}}$. Let $\ld_0 = \ld$, and let
\eq\label{eq:ldj}
\ld_j = s_{\af_{i_j}} \cdot \ld_{j-1}, \quad
\beta_j = s_{\af_{i_1}} \cds s_{\af_{i_{j-1}}}(\af_{i_j})
~\textup{ for }~j=1,\ldots, N.
\endeq
Note that for each $j$, the map 
$M(\ld_j)_\CC \rw M(\ld_{j-1})_\CC$ sending 
$v^+_{\ld_j}$ to 
$S(\af_{i_j}, \g{\ld_{j-1}+\rho, \af_{i_j}^\vee}) v^+_{\ld_{j-1}}$
is an inclusion, where
\[
\g{\ld_{j-1}+\rho, \af_{i_j}^\vee}
= 
\g{s_{\af_{i_{j-1}}} \cds s_{\af_{i_{1}}}(\ld+\rho), \af_{i_j}^\vee}
=
\g{\ld+\rho, \beta_j^\vee}.
\]
Therefore, the composition
\[
M(w\cdot\ld)_\CC \rw M(\ld_{N-1})_\CC \rw \cds \rw M(\ld_1)_\CC \rw M(\ld)_\CC
\]
is again an inclusion sending $v^+_{w\cdot \ld}$ to $S v^+_\ld$,
where $S$ is a product of generic Shapovalov elements, given by
\eq\label{eq:S}
S =
S(\af_{i_1}, \g{\ld+\rho, \beta_1^\vee}) 
\cds
S(\af_{i_N}, \g{\ld+\rho, \beta_N^\vee}) 
=
f_{\af_{i_1}}^{\g{\ld+\rho, \beta_1^\vee}}
\cds
f_{\af_{i_N}}^{\g{\ld+\rho, \beta_N^\vee}}
.
\endeq
\begin{lem}\label{lem:FSvk}
Assume that $\ld\in X^+, w\in W_{I\backslash\{\eta\}}$ and $\beta \in Q^+$. 
Then for each $k \in \ZZ_{\geq 0}$, 
\[
\{F_\pi S v^+_{\ld+k\varpi_\eta} ~|~  \textup{wt}(\pi) = w\cd\ld - (\ld-\beta)\}
\]
is a basis of the image of $M(w\cdot \ld+k\varpi_\eta )_\CC$ in $M(\ld+k\varpi_\eta)_\CC$,
where $S$ is as defined in \eqref{eq:S}.
\end{lem}
\proof
Fix a reduced expression $w = s_{\af_{i_N}}\cds s_{\af_{i_1}}$ and define $\ld_j, \beta_j$ as in \eqref{eq:ldj}. Since $w \in W_{I\backslash\{\eta\}}$, we have
 $w(\varpi_\eta) = \varpi_\eta$ 
and hence for each $k\in \ZZ$ and $j=1,\ldots, N$, 
\[
w\cd(\ld_j+k\varpi_\eta) = w\cd\ld_j + w(k\varpi_\eta) =  w\cd\ld_j + k\varpi_\eta.
\]
By definition of generic Shapovalov elements, the map 
$M(\ld_j+k\varpi_\eta)_\CC \rw M(\ld_{j-1}+k\varpi_\eta)_\CC$ sending 
$v^+_{\ld_j+k\varpi_\eta}$ to 
$S(\af_{i_j}, \g{\ld+\rho, \beta_j^\vee}) v^+_{\ld_{j-1}+k\varpi_\eta}$
is an inclusion for each $j$, which leads to another inclusion
$
M(w\cdot\ld+k\varpi_\eta)_\CC \rw  M(\ld+k\varpi_\eta)_\CC
$
sending $v^+_{w\cdot \ld+k\varpi_\eta}$ to $S v^+_\ld+k\varpi_\eta$, where $S$ is exactly the same element defined in \eqref{eq:S}. 
\endproof
\lem \label{lem:basis}
Assume that 
$\{B_i \overline{v}^+_\ld\}_{i\in \Ld}$ is a basis of $L(\ld)_{\CC,\ld-\beta}$ with $B_i \in U(\fn_\CC^-)$
for some $\ld\in X^+$ and 
$\beta=\sum m_i\af_i \in Q^+$.
If $\g{\ld+\rho,\af_\eta^\vee} > m_\eta$, then for each $k\in \ZZ_{\geq 0}$,
\[
\{B_i \overline{v}^+_{\ld+k\varpi_\eta}\}_{i\in \Ld}
\]
 is a basis of $L(\ld+k\varpi_\eta)_{\CC,\ld+k\varpi_\eta-\beta}$ satisfying the property \eqref{*} below:
for each $y\in U(\fn_\CC^-)_{-\beta}$ and $k\in \ZZ_{\geq 0}$, there are  $C^{(y)}_i \in \CC$ (independent of $k$) such that
\eq\label{*}
y \overline{v}^+_{\ld+k\varpi_\eta} = \sum_{i\in\Ld} C^{(y)}_i B_i \overline{v}^+_{\ld+k\varpi_\eta}.
\endeq
\endlem
\proof
By Lemma \ref{lem:wt} and Lemma \ref{lem:FSvk}, 
there is a basis $\{B'_j \overline{v}^+_\ld\}_{j\in \Ld'}$ of $L(\ld)_{\CC,\ld-\beta}$ with $B'_j \in U(\fn_\CC^-)$ satisfying Property \eqref{*},
 since $F_\pi S$ is independent of $k$.
Plugging in $y = \sum_{i\in \Ld} C_iB_i$ for arbitrary $C_i \in \CC$ to \eqref{*}, we obtain that 
\[
\sum_{i\in \Ld}C_i B_i \overline{v}^+_{\ld+k\varpi_\eta} 
= \sum_{j\in \Ld'}C'_j B'_j \overline{v}^+_{\ld+k\varpi_\eta},
\] 
where 
$C'_j = C^{(\sum_{i\in \Ld} C_iB_i)}_j$.

Note that $\sum_{j\in \Ld'}C'_j B'_j \overline{v}^+_{\ld+k\varpi_\eta}= 0$ if and only if all the $C'_j$'s are zero, which is equivalent to that all  $C_i$'s are zero. 
Therefore $\{B_i \overline{v}^+_{\ld+k\varpi_\eta}\}_{i\in \Ld}$ is a basis of $L(\ld+k\varpi_\eta)_{\CC,\ld+k\varpi_\eta-\beta}$ for each $k\in \ZZ_{\geq 0}$. 
Moreover, $\{B_i \overline{v}^+_{\ld}\}_{i\in \Ld}$ also has Property \eqref{*}.
\endproof
\subsection{Homomorphisms between Weyl Modules}
Here we fix $\ld,\mu = \ld - D\gamma \in X^+$ satisfying $\mu\uw_\gamma \ld$ for some $\eta$-good positive real root $\gamma=\sum g_i \af_i$ 
and $D\in \ZZ_{>0}$. 
Fix $Z_\eta = Z_\eta(\gamma,D)$ to be an $h_\eta$-avoiding integral Shapovalov element.
Fix a basis $\{B_a \overline{v}^+_\ld\}_{a\in \mc{B}_1}$ of $L(\ld)_{\ZZ,\mu}$
with $B_a \in U(\fn_\ZZ^-)_{-D\gamma}$,
and fix a basis $\{B'_b \overline{v}^+_\ld\}_{b\in \mc{B}_2}$ of $L(\ld)_{\ZZ, \mu+n\af_i}$ for some fixed $i\in I$ and $n\in\ZZ_{\geq 0}$, 
with $B'_b \in U(\fn_\ZZ^-)_{-D\gamma+n\af_i}$. 

We start with a corollary of Lemma \ref{lem:basis}.
\cor \label{cor:basis}
If $\g{\ld+\rho,\af_\eta^\vee} > Dg_\eta$, then for each $k \in \ZZ_{\geq 0}$ we have
\eq\label{2nd}
e_i^{(n)} Z_\eta \overline{v}^+_{\ld+ k\varpi_\eta} 
= \sum_{m=1}^{Dg_\eta} \sum_{b\in\mathcal{B}_2}
u_{b,m}(\ld) 
\ts{\g{\ld+k\varpi_\eta+\rho,\gamma^\vee}-D \choose m} 
B'_b \overline{v}_{\ld+k\varpi_\eta}^+,
\endeq
for some polynomials $u_{b,m} \in \CC[h_0,\ldots,\widehat{h}_\eta,\ldots,h_r]$.
\endcor
\proof
Applying Lemma \ref{lem:basis} with $\beta = D\gamma - n\af_i$ and $y = F_\pi$ with wt$(\pi) = D\gamma-n\af_i$,  we have, for each $k\in\ZZ_{\geq0}$,
\[
F_\pi \overline{v}^+_{\ld+ k\varpi_\eta} 
= \sum_{b\in \mathcal{B}_2} C_b^{(F_\pi)} B'_b \overline{v}^+_{\ld+ k\varpi_\eta},
\]
where $C_b^{(F_\pi)}$ is independent of $k$.
Let $u_m \in U(\fn_\ZZ^-) \otimes \CC[h_0, \ldots, \widehat{h}_\eta, \ldots, h_r]$ be as defined in \eqref{eq:PbeZ}. We have
\[
u_m = \sum_{\pi \in \mathcal{P}} F_\pi u_{m,\pi},
\]
for some $u_{m,\pi} \in \CC[h_0, \ldots, \widehat{h}_\eta, \ldots, h_r]$,
that is, $u_{m,\pi}(\ld+k\varpi_\eta) = u_{m,\pi}(\ld)$ for any $k\in\ZZ$. 
Thus, by Proposition \ref{prop:PbeS}, we have
\[
\ba{{ll}
\ds e_i^{(n)} Z_\eta \overline{v}^+_{\ld+ k\varpi_\eta} 
= \sum_{m=1}^{Dg_\eta} u_m 
\ts{ h_\gamma + \rho(h_\gamma) -D \choose m} 
\overline{v}^+_{\ld+ k\varpi_\eta} 
\\
\ds=
\sum_{m=1}^{Dg_\eta} 
\sum_{\pi \in \mathcal{P}}  
u_{m,\pi} (\ld)
\ts{ \g{\ld+k\varpi_\eta+ \rho, \gamma^\vee} -D \choose m} 
F_\pi
\overline{v}^+_{\ld+ k\varpi_\eta} 
\\
\ds=
\sum_{m=1}^{Dg_\eta} 
\sum_{b\in \mathcal{B}_2} 
\left(
\sum_{\pi \in \mathcal{P}}  
u_{m,\pi} (\ld)
C_b^{(F_\pi)} 
\right)
\ts{ \g{\ld+k\varpi_\eta+ \rho, \gamma^\vee} -D \choose m} 
B'_b \overline{v}^+_{\ld+ k\varpi_\eta}.
}
\]
Let $u_{m,b} = \sum_{\pi \in \mathcal{P}}  
u_{m,\pi}
C_b^{(F_\pi)} $  and we are done.
\endproof


\lem \label{lem:main}
If $\g{\ld+\rho,\af_\eta^\vee} > Dg_\eta$, then for each $k \in \ZZ_{\geq 0}$, we have
\[
e_i^{(n)} 
Z_\eta 
\overline{v}^+_{\ld+ k\varpi_\eta} 
= 
\sum_{m=1}^{Dg_\eta} 
\sum_{\substack{a\in \mathcal{B}_1\\b\in\mathcal{B}_2}}
C_a C_{a,b,m}(\ld) 
\ts{\g{\ld+ k\varpi_\eta+\rho,\gamma^\vee}-D \choose m} 
B'_b 
\overline{v}_{\ld+ k\varpi_\eta}^+,
\]
for some $C_a \in \CC$ and polynomials $C_{a,b,m} \in \CC[h_0,\ldots,\widehat{h}_\eta,\ldots,h_r]$,
which are independent of $k$.
\endlem
\proof
Now let $Z_\eta\overline{v}_\ld^+ = \sum_{a\in\mathcal{B}_1} C_a B_a \overline{v}^+_\ld$ with $C_a\in \CC$. 
Applying Lemma \ref{lem:basis} with $\beta = D\gamma - n\af_i$ and $y = e_i^{(n)} B_a$,  we have, for each $k\in\ZZ_{\geq0}$,
\[
e_i^{(n)} B_a \overline{v}^+_{\ld+ k\varpi_\eta} 
= \sum_{b\in \mathcal{B}_2} C_b^{(e_i^{(n)} B_a)} B'_b \overline{v}^+_{\ld+ k\varpi_\eta},
\]
where each $C_b^{(e_i^{(n)} B_a)}$ is independent of $k$. Similar to the proof of Corollary \ref{cor:basis}, one obtains that
\[
e_i^{(n)} B_a \overline{v}^+_{\ld+ k\varpi_\eta} 
= \sum_{m\geq0} \sum_{b\in\mathcal{B}_2}
 C_{a,b,m}(\ld) 
\ts{\g{\ld+ k\varpi_\eta+\rho,\gamma^\vee}-D \choose m} 
B'_b 
\overline{v}^+_{\ld+ k\varpi_\eta},
\]
for some $C_{a,b,m}\in \CC[h_0,\ldots,\widehat{h}_\eta,\ldots,h_r,d]$ that is independent of $k$.
Therefore
\eq
e_i^{(n)} Z_\eta \overline{v}^+_{\ld+ k\varpi_\eta} 
= \sum_{m\geq0} \sum_{\substack{a\in\mathcal{B}_1\\b\in\mathcal{B}_2}}
 C_aC_{a,b,m}(\ld) 
\ts{\g{\ld+ k\varpi_\eta+\rho,\gamma^\vee}-D \choose m} 
B'_b 
\overline{v}^+_{\ld+ k\varpi_\eta}.
\endeq\label{eq:1st}
Equating the coefficient of each $B'_b\overline{v}^+_{\ld+k\varpi_\eta}$ on \eqref{2nd} and \eqref{eq:1st}, one gets
\[
\ba{{l}
\ds\sum_{m\geq0} 
\sum_{a\in\mathcal{B}_1} 
C_aC_{a,b,m}(\ld)  
\ts{\g{\ld+k\varpi_\eta+\rho,\gamma^\vee}-D \choose m}
\ds=
\sum_{m=1}^{Dg_\eta} 
u_{b,m}(\ld) 
\ts{\g{\ld+k\varpi_\eta+\rho,\gamma^\vee}-D \choose m}
\in\CC[k].
}
\]
Since the right hand side has finite degree in $k$, 
we may equate the coefficients of powers of $k$, and hence
\[
\sum_{a\in\mathcal{B}_1} C_aC_{a,b,m}(\ld) = 
\bc{
u_{b,m}(\ld)
&\tif 1\leq m \leq Dg_\eta,
\\
0&\tif m =0 \textup{ or } m> Dg_\eta.
}
\]
\endproof
Now we can prove our main theorem.
\thm \label{thm:main}
Assume that $\mu \uw^e_\gamma \ld$ and $\gamma=\sum g_i \af_i \in \Phi_{\text{re}}^+$ is $\eta$-good.
If $Dg_\eta < \g{\ld+\rho,\af_\eta^\vee}$  and $Dg_\eta < p^e$
then there exists a nonzero homomorphism $V(\mu)\rw V(\ld)$  between Weyl modules.
\endthm
\proof
Let $p^g$ be the  highest $p$-power dividing $Z_\eta\overline{v}_\ld^+ = \sum_{a \in \mathcal{B}_1} C_a B_a \overline{v}_\ld^+$ in $L(\ld)_\ZZ$,
that is, $C_a \in p^g\ZZ$ for each $a$.
We shall show that the map sending $\overline{v}_\mu^+\otimes 1_\KK$ to $\frac{Z_\eta\overline{v}_\ld^+}{p^g}\otimes 1_\KK$ is a homomorphism.
Applying Lemma \ref{lem:main}, we have
\[
e_i^{(n)} \frac{Z_\eta \overline{v}^+_{\ld}}{p^g} 
= \sum_{m=1}^{Dg_\eta} 
\sum_{\substack{a\in\mathcal{B}_1 \\ b\in\mathcal{B}_2}}
\frac{C_a}{p^g} C_{a,b,m}(\ld) 
\ts{\g{\ld+\rho,\gamma^\vee}-D \choose m} 
B'_b \overline{v}_{\ld}^+.
\]
Since $Dg_\eta < p^e$, each $m <p^e$ and hence each ${\g{\ld+\rho,\gamma^\vee}-D \choose m} = {Mp^e \choose m} \in p\ZZ$. 
Therefore $e_i^{(n)} \frac{Z_\eta \overline{v}^+_{\ld}}{p^g} \otimes 1_\KK = \overline{0}$ in $V(\ld)$.
\endproof
\exa
In Example \ref{exa4} we have $\gamma = 2\af_0 + \af_1 = 2\delta - \af_1 = 2\varpi_0 - 2\varpi_1 + 2\delta$ is $1$-good for any prime with $g_1 =1$ and
\[
\ba{{l}
Z_1 = f_1f_0^2 - f_{\delta,1} f_0 (h_0-1) + f_\gamma (h_0^2-h_0)
\\
= f_1 f_0^{(2)} (2h_0^2) - f_0f_1f_0 (2h_0^2 -h_0 -1) + f_0^{(2)}f_1 (2h_0^2 -2h_0).
}
\]
For $\ld= 2\varpi_0+\varpi_1$, then $\mu = \ld - D\gamma \in X^+$ implies $D=1$ and $\mu = 3\varpi_1 - 2\delta$.
 Note that $\g{\ld+\rho,\gamma^\vee} = 8 = Mp^e + D$.
 Hence $p=7$ is the only possible prime with $M=e=1$. 
Also, the assumptions are satisfied since $Dg_1 = 1 <\g{\ld+\rho,\af_1^\vee} = 2$ and $Dg_1 = 1 < p^e = 7$.
Note that
\[
Z_1 \overline{v}^+_\ld 
=\left( 8 f_1 f_0^{(2)}  - 5 f_0f_1f_0  + 4 f_0^{(2)}f_1 \right)\overline{v}^+_\ld,
\]
so $g=0$
and therefore the map 
$
V(3\varpi_1-2\delta) \rw V(2\varpi_0+\varpi_1)
$ 
sending $\overline{v}_\mu^+\otimes 1_\KK$ to $Z_1\overline{v}_\ld^+\otimes 1_\KK$ is a nonzero homomorphism.
\endexa
\section{A conjectural strong linkage principle}  \label{sec:red Weyl}
\subsection{Formulation of a conjecture}
Having in mind the strong linkage principles for classical cases, 
we formulate a conjecture (joint with W. Wang) on the strong linkage principle for the modular representation of affine Lie algebras.
\conj[Strong linkage principle]\label{conj}
Let $\mu,\ld\in X^+$, and let $L(\mu)$ be the unique irreducible highest weight $U_\KK$-module with highest weight $\mu$.
If $L(\mu)$ is a composition factor of the Weyl module $V(\ld)$, 
then $\mu =\ld$ or $\mu = \mu_0 \uw \cds \uw \mu_N = \ld$ for some $N$ and $\mu_i \in X$ where
\eq\label{SL:aff mod}
x \uw y \Lrw \text{ there exists } 
\bc{
	n\in \ZZ_{>0},\\ 
	m\in \ZZ,\\ 
	\beta \in \Phi^+
} 
\text{such that } 
\bc{
	y-x= (n-mp)\beta \in Q^+,\\ 
	n(\beta,\beta) = 2(y+\rho, \beta) \in \ZZ/p\ZZ.
}
\endeq
\endconj
\rmk
If we interpret $p=0$, condition \eqref{SL:aff mod} coincides with Condition \eqref{SL:aff ord}.
For finite types, Condition \eqref{SL:aff mod} describes exactly the strong linkage for modular finite case below:
\eq
x \uw y \Lrw \text{ there exists } 
\bc{
	m\in \ZZ,\\ 
	\beta \in \Phi^+
} 
\text{such that } 
	x= s_{\beta,mp}\cd y < y.
\endeq
\endrmk
\rmk
The formulation of Conjecture \ref{conj} also makes sense for all symmetrizable Kac-Moody algebras.
\endrmk
\subsection{Candidates for reducible Weyl modules}
Now we focus on the reducibility problem for Weyl modules.
Note that we do not know the irreducibility for any Weyl module, but our main theorem is able to detect if a Weyl module is reducible.

The candidates for the high weights of reducible Weyl modules are those dominant integral weights which are $\gamma$-mirrored by another dominant integral weight. 
Note that if $\ld,\mu\in X^+$ such that $\mu\uw_\gamma \ld$ for some positive real root $\gamma$,
there is a unique weight $\ld+t\delta \in X^+$ such that $\mu+t\delta \uw_\gamma \ld+t\delta$ and $\xi = 0$ if we express $\ld+t\delta = \sum_i \xi_i \varpi_i + \xi\delta$. 
Without loss of generality, we only need to consider the set defined below:
 \[
Y^+=\left\{ \ld\in \sum_i\ZZ_{\geq 0}\varpi_i~\middle|~  \e \mu\in X^+ \text{ s.t. } \mu \uw_\gamma\ld \fsm \gamma\in\Phi_{\text{re}}^+ \right\}.
\] 
\lem\label{lem:phi}
For an arbitrary affine type, 
the set of positive real roots can be described as the following.
\[
\Phi^+_{\text{re}} = \bigcup_{i=1}^3 \left\{\gamma_0 + t\delta ~\middle|~ \gamma_0 \in \Phi_i^+, t\in i\ZZ_{\geq 0} \right\},
\]
where $\Phi_i^+$ are subsets of $\Phi^+$ for $i=1,2,3$.

In particular, $\Phi_3^+ = \varnothing$ unless for type $\widetilde{D}_4^{(3)}$.
$\Phi_2^+ = \varnothing$ unless for type $\widetilde{A}_r^{(2)}, \widetilde{D}_{r}^{(2)}$ and $\widetilde{E}_6^{(2)}$.
For untwisted types we have $\Phi^+_1 =\{\af, \delta-\alpha ~|~ \af\in \Phi^+_0\}$.
\endlem
\proof
This is done by a case-by-case analysis on the data of root systems.
\endproof
\prop\label{prop:Y+}
Let $\ld \in X^+$ be a dominant integral weight of fixed level $\ell$, and let $\gamma_0\in \Phi_i^+$ for some $i$ as defined in Lemma \ref{lem:phi}. Let $C_i = \frac{2i}{(\gamma_0,\gamma_0)}(\ell+h^\vee)$.
If $\ld - \gcd(C_i,p)\gamma_0\in X^+$ then $\ld \in Y^+$.
\endprop
\proof
For any positive real root of the form $\gamma = \gamma_0+it\delta$ we have $(\gamma,\gamma) = (\gamma_0,\gamma_0)$, 
and hence for any $\ld\in X^+$ we have $\g{\ld+\rho,\gamma^\vee} = \g{\ld+\rho,\gamma_0^\vee} + C_i t$.
Therefore for $m,e \in\ZZ_{>0}$ we have
\[
s_{\gamma,Mp^e} \cd\ld = s_{\gamma_0,C_it+Mp^e}\cd\ld - (\g{\ld+\rho,\gamma^\vee}-Mp^e)t\delta.
\]

Therefore, we can always choose $t,M \in \ZZ_{>0}$ large enough such that 
\[
s_{\gamma_0+it\delta,Mp^e} \cd\ld = \ld-\gcd(C_i,p^e)\gamma_0 - (\g{\ld+\rho,\gamma^\vee}-Mp^e)t\delta.
\]

By assumption $s_{\gamma_0+t\delta,Mp^e} \cd\ld \in X^+$ and hence $\ld\in Y^+$.
\endproof~\\

Using Proposition \ref{prop:Y+}, one can describe for each type and for each prime $p$ the lowest level weights in $Y^+$.
For most types, every weight in $Y^+$ has level $\geq 2$. $Y^+$ contains a level one weight if and only if the following conditions hold:
\enu
\item $p$ is a odd prime that does not  divide $h^\vee+1$.
\item $\fg_\CC$ is of type 	$\widetilde{B}_r$, $\widetilde{C}_r$, $\widetilde{F}_4$ or $\widetilde{G}_2$.
\endenu{
In this case, the level one weights are given by Table 1 below.

\begin{table}[!h]
\label{table:Y+}
\caption{The complete list of weights of level one in $Y^+$.}
\[
\ba{{|c|c|c|c|c|cccccccccccccccccccccc}
\hline&&&&\\
\text{Type}
	& \widetilde{B}_r; r\geq 3 
	& \widetilde{C}_r; r\geq 2 
	& \widetilde{F}_4
	& \widetilde{G}_2
\\
\hline&&&&\\
\text{weight}
	& \varpi_0, \varpi_1
	& \varpi_0, \varpi_r
	& \varpi_0
	& \varpi_0, \varpi_2
\\
\hline
}
\]
\end{table}

This shows that for (possibly twisted) affine ADE types, 
there is no dominant integral weights that are $\gamma$-mirrored to $\varpi_0$ for any $\gamma \in \Phi_{\text{re}}^+$. 
One might expect that the lack of $\gamma$-mirrored weights would imply that the basic representation is irreducible.
However, the table in \cite[\S1]{BK02} shows that for each (possibly twisted) affine ADE type, 
there exists a reducible basic representation for some prime $p\leq h$. 
\endrmk
\subsection{Quasi-simple weights}\label{defn:qs}
Our main theorem shows that the corresponding Weyl module of a weight in $Y^+$ is reducible if the restrictions  in Theorem \ref{thm:main} are satisfied. 
Here we shall demonstrate that these restrictions are actually quite mild.

\defn\label{defn:exc}
A weight $\ld\in Y^+$ is called \textit{quasi-simple} if Theorem \ref{thm:main} does not apply for any $\mu\in X^+$ that is $\gamma$-mirrored to $\ld$.
\enddefn

Note that every weight of level one described in Remark \ref{table:Y+} is quasi-simple. 
Unfortunately, it seems that there is not an obvious pattern for the quasi-simple weights in general.

\subsection{Type $\widetilde{A}_1$}
Assume that $\gamma = \af_0+t\delta$ with $t\geq 0$.
Let $(\xi_0,\xi_1)$ be the shorthand notation of the weight $\xi_0\varpi_0+\xi_1\varpi_1$.
Our main theorem has the following corollary.
\cor\label{cor:A1}
For type $\widetilde{A}_1$ and a fixed prime $p$, a weight $\ld =( \xi_0, \ell-\xi_0)$ of level $\ell$ lies in $Y^+$ if and only if there are integers
$M,D,e\in\ZZ_{>0}$ and $t\in \ZZ_{\geq0}$ satisfying the following conditions.
\enu
\item $(\ell+2)t+\xi_0+1=Mp^e+D$.
\item $2D \leq \xi_0 \leq \ell$.
\endenu

In particular, if $\ld\in Y^+$, the corresponding Weyl module $V(\ld)$ is reducible if either Condition (i) or (ii) below holds:
\enu 
\item[(i)] $\gcd(t+1,p)=1$ and $D(t+1) < \min(\xi_0+1,p^e)$.
\item[(ii)] $\gcd(t,p)=1$ and $Dt < \min(\ell-\xi_0+1,p^e)$.
\endenu
\endcor
\proof
This is a reformulation of Theorem \ref{thm:main} in type $\widetilde{A}_1$.
\endproof

For a prime $p$, we denote by $\ell\ell(p)$ the lowest level of weight $\ld$ such that $V(\ld)$ is a reducible Weyl module arising from our main theorem. 
Note that there could be reducible Weyl modules that cannot be detected by our main Theorem with level lower than $\ell\ell(p)$.
Table 2 below describes the first few $\ell\ell(p)$ and the possible high weights $\ld = (\xi_0,\ell\ell(p)-\xi_0)$ of reducible Weyl modules detected by Corollary \ref{cor:A1}. One can observe that $\ell\ell(p)$ grows much slower than $p$.
\begin{table}[!h]
\caption{The lowest level $\ell\ell(p)$ and possible $\xi_0$.}
\[
\ba{{c|cccccccccccccccccc}
p&2&3&5&7&11&13&17&19\\
\hline
\ell\ell(p)&2&2&4&3&4&3&5&5\\
\hline
\text{possible }\xi_0&\{0,2\}&\{0,2\}&\{0,1,3,4\}&\{1,2\}&\{0,4\}&\{0,3\}&\{2,3\}&\{0,5\}\\
\\
p&23&29&31&37&41&43&47&53\\
\hline
\ell\ell(p)&5&6&7&6&7&7&7&8\\
\hline
\text{possible }\xi_0&\{2,3\}&\{1,5\}&\{3,4\}&\{1,5\}&\{2,5\}&\{0,7\}&\{2,5\}&\{3,5\}\\
}
\]
\end{table}

The quasi-simple weights for $p = 2,3,$ and $5$ up to level 150 are given in Table 3 below. 
In other words, any weight in $Y^+$ that is not in the table corresponds to a reducible Weyl module. 
This supports the hypothesis that almost every Weyl module is reducible.

\begin{table}[ht!]
\caption{Quasi-simple weights of level $<150$ in type $\widetilde{A}_1$.}
\[
\ba{{cccc}
p=2 &p=3&p=5\\
\bt{{| c| l| c| cccccccc}
\hline
$\ell$& $(\xi_0, \ell-\xi_0)$ \\
\hline
3&(0,3), (3,0) \\
8&(1,7), (7,1) \\
18&(3,15), (15,3) \\
38&(7,31), (31,7) \\
78&(15,63), (63,15) \\
\hline}
&
\bt{{| c| l| c| cccccccc}
\hline
$\ell$&$(\xi_0, \ell-\xi_0)$\\
\hline
3&(1,2), (2,1) \\
6&(1,5), (5,1) \\
13&(5,8), (8,5) \\
22&(5,17), (17,5) \\
43&(17,26), (26,17) \\
70&(17,53), (53,17) \\
\hline
}
&
\bt{{| c| l| c| cccccccc}
\hline
$\ell$& $(\xi_0, \ell-\xi_0)$\\
\hline
2&(0,2), (2,0) \\
4&(2,2) \\
6&(3,3) \\
18&(4,14), (14,4) \\
28&(14,14) \\
38&(19,19) \\
98&(24,74), (74,24) \\
148&(74,74) \\
\hline}
}
\]
\end{table}
\subsection{Type $\widetilde{A}_r$ : $r\geq 2$}
Assume that $\gamma = \af_0+t\delta$ with $t\geq 0$.
Let $(\xi_0,\xi_1, \ldots , \xi_r)$ be the shorthand notation of the weight $\sum_{i=0}^r \xi_i\varpi_i$.
Our main theorem leads to the following corollary.
\cor\label{cor:Ar}
For type $\widetilde{A}_r, r\geq 2$ and a fixed prime $p$,
 a weight $\ld =(\xi_0, \xi_1, \ldots,$ $\xi_{r-1}, \ell-(\xi_0+\ldots + \xi_{r-1}))$ of level $\ell$ lies in $Y^+$ if and only if there are integers
$M, D,e\in\ZZ_{>0}$ and $t\in \ZZ_{\geq0}$ satisfying the following conditions.
\enu
\item $(\ell+h^\vee)t+\xi_0+1=Mp^e+D$.
\item $2D \leq \xi_0 \leq \ell - (\xi_1+\ldots + \xi_{r-1})$.
\endenu

In particular, if $\ld\in Y^+$, the corresponding Weyl module $V(\ld)$ is reducible if  either Condition (i),(ii) or (iii) below holds:
\enu 
\item[(i)] $\gcd(t+1,p)=1$ and $D(t+1) < \min(\xi_0+1,p^e)$.
\item[(ii)] $\gcd(t,p)=1$ and $Dt < \min(\xi_i+1,p^e)$ for some $i=1,\ldots, r-1$.
\item[(iii)] $\gcd(t,p)=1$ and $Dt < \min(\ell-(\xi_0+\ldots+\xi_{r-1})+1,p^e)$.
\endenu
\endcor
\rmk\label{rmk:upper bdd}
We observe that,  for a fixed type, every weight of level $\ell$ seems to be quasi-simple if $p \gg \ell$. 
%
In particular, for type $\widetilde{A}_r$ every weight in $Y^+$ of level $\ell$ is quasi-simple if $p > \ell(\ell+h^\vee)-h^\vee$.
\endrmk

\end{document}